\documentclass[final,11pt]{elsarticle}
\usepackage[latin1]{inputenc}
\usepackage{amsmath,bm}
\usepackage{amsthm}
\usepackage{amsfonts}
\usepackage{color}
\usepackage[dvipsnames]{xcolor}
\usepackage{amssymb}
\usepackage{graphicx}
\usepackage{comment}
\usepackage{comment,booktabs,multirow}
\usepackage{gensymb}
\usepackage{fullpage}

\usepackage{caption,multirow,booktabs}
\usepackage{xfrac}
\usepackage{subcaption}
\usepackage[linesnumbered,ruled,vlined]{algorithm2e}

\usepackage{empheq}
\usepackage[top=1in, bottom=1in, left=1in, right=1in]{geometry}

\usepackage[pdftex, unicode=true,
colorlinks,
linkcolor=red,
anchorcolor=blue,
citecolor=green
]{hyperref}

\newcommand{\ahalf}{\sfrac{1}{2}}
\newtheorem{thm}{Theorem}[section]

\newtheorem{remark}[thm]{Remark}
\newtheorem{proposition}{Proposition}

\newcommand{\lb}{\left (}
\newcommand{\rb}{\right )}

\newcommand{\bu}{\bm u}

\newcommand{\bV}{\bm V}

\newcommand{\bk}{\bm k}

\newcommand{\bF}{\bm F}

\makeatletter
\newcounter{manualsubequation}
\renewcommand{\themanualsubequation}{\alph{manualsubequation}}
\newcommand{\startsubequation}{%
	\setcounter{manualsubequation}{0}%
	\refstepcounter{equation}\ltx@label{manualsubeq\theequation}%
	\xdef\labelfor@subeq{manualsubeq\theequation}%
}
\newcommand{\tagsubequation}{%
	\stepcounter{manualsubequation}%
	\tag{\ref{\labelfor@subeq}\themanualsubequation}%
}
\let\subequationlabel\ltx@label
\makeatother

\makeatletter
\def\ps@pprintTitle{%
	\let\@oddhead\@empty
	\let\@evenhead\@empty
	\def\@oddfoot{}%
	\let\@evenfoot\@oddfoot}
\makeatother

\SetKwInput{KwInput}{Input}                
\SetKwInput{KwOutput}{Output}              

\DeclareCaptionLabelFormat{algnonumber}{Algorithm}
\begin{document}
	\begin{frontmatter}
		\title{High-Order Multirate Explicit Time-Stepping Schemes \\for the Baroclinic-Barotropic Split Dynamics in Primitive Equations}
		\author[usc]{Rihui Lan}
		\ead{rlan@mailbox.sc.edu}
		\author[usc]{Lili Ju\corref{cor}}
		\ead{ju@math.sc.edu}
		\author[usc]{Zhu Wang}
		\ead{wangzhu@math.sc.edu}
		\author[fsu]{Max Gunzburger}
		\ead{mgunzburger@fsu.edu}
		\author[lanl]{Philip Jones}
		\ead{pwjones@lanl.gov}
		\address[usc]{Department of Mathematics, University of South Carolina, Columbia, SC 29208, USA}
		\address[fsu]{Department of Scientific Computing, Florida State University, Tallahassee, FL 32306, USA}
		\address[lanl]{Theoretical Division, Los Alamos National Laboratory, Los Alamos, NM 87545, USA}
		\cortext[cor]{Corresponding author}
		\begin{abstract} 
		
In order to treat the multiple time scales of ocean dynamics in an efficient manner, the baroclinic-barotropic splitting technique has been widely used for solving the primitive equations for ocean modeling. Based on the framework of strong stability-preserving Runge-Kutta  approach, we propose two high-order  multirate  explicit 
time-stepping schemes (SSPRK2-SE and SSPRK3-SE) for the resulting split system in this paper. 
The  proposed schemes allow for a large time step  to be used for  the three-dimensional  baroclinic (slow) mode and a small time step for the two-dimensional barotropic (fast)  mode,  in which each of the two mode solves just need to satisfy their respective CFL conditions for numerical stability. Specifically, at each time step, the  baroclinic velocity is first computed  
by advancing the baroclinic mode and fluid thickness of the system with the large time-step   \textcolor{black}{and the assistance of  some intermediate approximations of the baroctropic mode obtained by substepping with the small  time step};
then the barotropic velocity is corrected by using  the small time step to re-advance the barotropic mode under an improved barotropic forcing produced by interpolation of the forcing terms from the preceding baroclinic mode solves; lastly, the fluid thickness is updated by coupling the baroclinic and barotropic velocities. Additionally, numerical inconsistencies on the discretized sea surface height caused by the mode splitting  are  relieved  via a reconciliation process with carefully calculated flux deficits. 
Two benchmark tests from the ``MPAS-Ocean" platform are  carried out to numerically demonstrate  the performance and parallel scalability of the proposed SSPRK-SE schemes.
		\end{abstract}
		
		\begin{keyword}
			Primitive equations,  baroclinic-barotropic splitting, explicit time-stepping, multirate, strong stability preserving Runge-Kutta, SSH reconciliation
		\end{keyword}
	\end{frontmatter}
	
	\section{Introduction}

The importance of the ocean to our everyday lives, now and more so in the future, is undeniable because it has a huge impact on the climate. To better understand and forecast the ocean and its effect on climate, a number of computational models, based on fundamental physics laws and the properties of geophysical flows, have been developed and used in practical applications. Among them are the primitive equations \cite{chen2003unstructured, smagorinsky1963general, vallis2006} which are a simplification of the Navier-Stokes equations. Those equations couple tracers such as temperature, salinity, and chemicals to the fluid velocity, depth (or layer thickness), and pressure. The shallow-water equations \cite{kinnmark2012shallow, vallis2006,  vreugdenhil1994numerical} which are a further simplification, are often employed in some specific circumstances such as for flows in rivers and coastal areas. Due to the high horizontal-to-vertical aspect ratio of the ocean,  models are generally derived by modeling the fluid as a single layer or as a stack of immiscible layers, each having a uniform fluid density. Consequently, layered models are ideal for modeling stratified fluid flows and perform well at portraying vertical profiles. For additional details, we refer to \cite{cushman2011introduction, vallis2006} and the references cited therein.

Numerical simulations of the layered models are still challenging due to the massive computational complexity caused by integrating the large-scale dynamical system. Because external and internal gravity waves and the Earth's rotation are intertwined, ocean dynamics often involves several different time scales. Explicit time-stepping schemes have been popularly applied in many numerical ocean models due to their natural parallelism and ease of implementation.  However, in order to achieve stable numerical simulations, such schemes have to use quite small time-step sizes imposed by the Courant-Friedrichs-Lewy (CFL) condition for the fastest time scales of the dynamics. This renders numerical simulations to be  computationally expensive, especially for long-term predictions. Therefore, it is natural to separate the modes of distinct characteristic time scales and advance them using different time-step sizes. To accomplish this goal, one approach is to develop explicit local time-stepping (LTS) schemes \cite{diaz2009energy, hoang2020high, HOANG2019152, trahan2012local}. Because CFL conditions vary greatly over the whole ocean domain, LTS methods apply spatially-dependent time-step sizes on different subdomains to achieve better efficiency compared to that obtained using globally uniform time-stepping methods. The other approach is to design multirate time-stepping schemes \cite{bleck1990wind, higdon1999implementation,higdon2002two,higdon2005two} that are based on splitting the ocean dynamics.  A large time-step size is used for advancing the slow dynamics mode and a small step size (i.e., substepping) for solving for the fast dynamics mode, so that each of the mode solves only needs to satisfy their respective CFL conditions for numerical stability.  In this paper, we focus on a novel instance of the latter approach and develop accurate multirate explicit time-stepping schemes  for the primitive equations.

Among the fast components of the primitive equations, the barotropic mode is the fastest of the entire spectrum of inertial-gravity waves. It has been a standard practice in ocean modeling to split that mode from the rest of the waves that make up the baroclinic mode. Based on the baroclinic-barotropic splittings, split-explicit (SE) time integration schemes have been developed in \cite{bleck1990wind, bryan1997numerical, dukowicz1994implicit,  higdon2002two, higdon2005two, higdon1997barotropic, killworth1991development, shchepetkin2005regional}. 
In those approaches, the barotropic velocity is obtained via the mass-weighted vertical averaging and the baroclinic velocity consists of the difference between the original velocity and the barotropic one. The whole dynamic system is then divided into a two-dimensional barotropic subsystem for the fast external gravity wave and a three-dimensional baroclinic subsystem for the slow baroclinic mode, whose time scales could be up to a hundred times different. 	
This splitting idea is used in \cite{bleck1990wind} to numerically predict the Atlantic ocean under the isopycnic coordinate system. However, in \cite{higdon1996stability} it was discovered that the scheme developed in \cite{bleck1990wind}  was not stable due to a certain inexactness in the splitting. Subsequently, a more stable two-level time-stepping method was proposed in \cite{higdon1997barotropic}. Furthermore, the splitting method was generalized to multi-layer ocean models in \cite{higdon1999implementation,higdon2002two,higdon2005two}. In particular, the barotropic mode in theses methods is obtained by averaging the momentum equation. For simplicity, only the linear Coriolis term and the gradient of the mass-weighted Montgomery potential are kept; all the other terms are then gathered in a special forcing term $\overline {\bm G}$ (usually called ``barotropic forcing''). Without accessing its (complicated) explicit formula, the value of $\overline {\bm G}$ is calculated efficiently from the baroclinic subsystem by using the mass-weighted free condition, so that  $\overline {\bm G}$  can bridge the baroclinic and barotropic modes. 
Furthermore, in \cite{higdon2005two}, a forward-backward substepping solver for the fast barotropic subsystem with $\overline {\bm G}$ calculated by the baroclinic mode. 
This two-level SE time-stepping approach is numerically very stable and has been implemented in  ``MPAS-Ocean'' \cite{ringler2013multi}, a numerical ocean model developed at  Los Alamos National Laboratory and its collaborating institutions  for the simulation of the ocean system across scales on staggered C-grids. This model was developed for use in the Energy Exascale Earth System Model (E3SM; \cite{E3SMv1}) as well as  some ocean-only applications \cite{Petersen_2019, reckinger2015study, wolfram2015diagnosing, woodring2015situ}. Implementation details about the MPAS-SE are given in \cite{petersen2018mpas}. The MPAS-SE package is based on the numerical schemes developed in \cite{ringler2013multi, ringler2010unified, thuburn2009numerical} where it is shown how the MPAS-SE scheme runs very stably and improves the efficiency of ocean simulations. However, the MPAS-SE scheme  is only of first-order accuracy in time as verified by our experiments in this paper. 
This is partly because the MPAS-SE only adopts the newly obtained $\overline {\bm G}$ and keeps it frozen during the barotropic time substepping while the Coriolis and the pressure-gradient terms are evolving. 
 
For the efficiently handling of the multiple time scales while simultaneously achieving higher-order accuracy for a given system, some variants of the  Runge-Kutta methods have been developed. In \cite{gunther2001multirate}, the multirate partitioned Runge-Kutta (MPRK) method was proposed and order conditions were derived based on the P-series \cite{hairer1993solving}. Recently, in \cite{sandu2015generalized, gunther2016multirate}, a generalized additive Runge-Kutta (GARK) method was designed which allows the evolution of the fast and slow modes using different numbers of stages and distinct sizes of time steps. However, neither the MPRK nor the GARK are directly applicable to the baroclinic-barotropic split system of the primitive equations for ocean modeling. One of the main difficulties is  that one would have to compute the barotropic forcing $\overline {\bm G}$ for each substep that advances the barotropic subsystem at all stages of MPRK or GARK, but the explicit form of $\overline {\bm G}$ is difficult to access in this case.
Some other classic high-order  time-stepping schemes have also been adopted in existing numerical ocean models \cite{griffies2000developments}, including   Leap-Frog, Adams-Bashforth, and forward-backward algorithms. In addition, the  TVD (total variation diminishing) \cite{shu1988total, leveque1992numerical,gottlieb1998total} and SSPRK (strong stability preserving Runge-Kutta) \cite{gottlieb2001strong, gottlieb2011strong} methods  have attracted some attention. In \cite{pietrzak1998use}, TVD limiters were considered for the advection of temperature and salinity in ocean modeling. In \cite{shi2012high}, SSPRK method was applied for the Boussinesq model to simulate wave shoaling, breaking, wave run-up, and wave-averaged nearshore circulation. In \cite{weller2013runge}, the SSPRK method was used with horizontally explicit, vertically implicit schemes to solve the compressible Boussinesq equations. 

Based on the SSPRK framework,  we propose in this paper two high-order  multirate explicit time-stepping schemes (SSPRK2-SE and SSPRK3-SE) for solving the baroclinic-barotropic split system of primitive equations, and each of the mode solves only needs to satisfy their respective CFL conditions.
Our approaches
		 are in some extent similar to the so-called compound-fast Runge-Kutta method \cite{verhoeven2006general,savcenco2007multirate,roberts2021implicit}.
 The SSPRK schemes are well known to achieve high-order temporal accuracy by utilizing a convex combination of forward-Euler steppings (stages). Through this strategy, we are able to compute the barotropic forcing  $\overline {\bm G}$ efficiently using the  approach of \cite{higdon2002two} at each stage of the SSPRK schemes. In particular, at each time step, the  baroclinic velocity is first computed  
by using the SSPRK scheme  to advance the baroclinic mode and fluid thickness of the system with the large time step   
\textcolor{black}{and the assistance of  some intermediate approximations of the baroctropic mode obtained by substepping with the small  time step;} 
Then the barotropic velocity is corrected by using the same SSPRK scheme with the small time step to re-advance the barotropic subsystem with a improved barotropic forcing produced by interpolation of the forcing terms from the preceding forward-Euler baroclinic mode solves. Lastly, the fluid thickness and the sea surface height (SSH) perturbation are updated by coupling the baroclinic and barotropic velocities.
	In the baroclinic-barotropic split system, either the barotropic mode or the layer thickness can be used to compute the SSH perturbation. However,  \textcolor{black}{they usually do not produce the same values at each stage of the SSPRK-SE schemes in the (time and/or space) discrete level},  which still could cause model inconsistency errors and  numerical stability issues  \cite{bleck1990wind,hallberg1997stable,higdon2005two,hallberg2009reconciling} to the whole discrete system. 
	There exist two typical ways to resolve this issue: the layer dilation \cite{bleck1990wind} and the flux-form reconciliation \cite{higdon2005two,hallberg2009reconciling}. The former is to uniformly dilate or compresses all the layers by applying a scaling factor on the related layer thickness equations, which reconciles the perturbation, but cannot conserve the mass of each layer  since the layer thickness equations are not cast in  flux-divergence form. The latter forces the sum of the layer fluxes to agree with the barotropic accumulated fluxes and the resulting equations are still in  flux-divergence form, thus it conserves the mass and is more preferred in practice.  For the proposed SSPRK-SE schemes, we design a new reconciliation strategy which is similar to the one developed in \cite{higdon2005two} but  uses  explicit transport velocity adjustments to better fit the  SSPRK framework.  

The rest of paper is structured as follows. In Section \ref{sect:dynamics}, we review the primitive equations and the corresponding baroclinic-barotropic dynamics splitting. In Section \ref{sect:ssprk}, we present the two multirate explicit time-stepping schemes, SSPRK2-SE and SSPRK3-SE, for the baroclinic-barotropic split system, and their error analyses in time are then discussed in Section \ref{sec:err}. 
In  Section \ref{sect:experiments}, numerical experiments based on  benchmark test cases from the ``MPAS-Ocean''  platform are carried out to demonstrate the performance and parallel scalability of the proposed schemes.  Finally, some concluding remarks are given in Section \ref{sect:conclusion}.
	
\section{The primitive equations and the baroclinic-barotropic dynamics splitting}\label{sect:dynamics}

We consider the primitive equations which describe  the incompressible Boussinesq equations in hydrostatic balance \cite{ringler2013multi,petersen2013mpas,petersen2015evaluation}. The  model consists of the following equations \cite{ringler2013multi}:
{\color{black}
	\begin{itemize}
		\item Thickness equation:
		\begin{equation}
		\frac{\partial h}{\partial t}+\nabla\cdot(h\bu)+\frac{\partial}{\partial z}(h w)=0;\label{layerthickness:single}
		\end{equation}
		\item Momentum equation:
		\begin{equation}
		\frac{\partial \bu}{\partial t}+\frac{1}{2}\nabla |\bu|^2+(\vec\bk\cdot\nabla\times\bu)\bu^{\perp}+f\bu^{\perp}+w\frac{\partial\bu}{\partial z}=-\frac{1}{\rho_0}\nabla p+\nu_h\nabla^2\bu+\frac{\partial}{\partial z}(\nu_v\frac{\partial\bu}{\partial z});\label{layer:momentum:single}
		\end{equation}
		\item Tracer equations:
		\begin{equation}
		\frac{\partial h\varphi}{\partial t}+\nabla\cdot(h\varphi\bu)+\frac{\partial}{\partial z}(h\varphi w)=\nabla\cdot(h\kappa_h\nabla\varphi)+h\frac{\partial}{\partial z}(\kappa_v\frac{\partial \varphi}{\partial z});\label{tracer:single}
		\end{equation}
		\item Hydrostatic condition:
		\begin{equation}
		p=p^s(x,y)+\int_z^{z^s}\rho g\,{\rm d}z';\label{hydrocon:single}
		\end{equation}
		\item Equation of state:
		\begin{equation}
		\rho=f_{\text{eos}}(\Theta,S,p),\label{eqstate:single}
		\end{equation}
\end{itemize}}
\noindent where the definitions of variables are listed in Table \ref{variable:def} and $\bu^{\perp}=\vec\bk\times\bu$.

\begin{table}[ht!]\small
	\centering
	\begin{tabular}{cl}
		\hline 
		Variables& Definition\\
		\hline
		$h$& Fluid thickness ({\color{black} total depth of a fluid column})\\
		$\bu$, $w$& Horizontal and vertical velocity\\
		$p$, $p^s(x,y)$& Pressure and the top surface pressure\\
		$\Theta$& Potential temperature\\
		$S$& Salinity\\
		$\rho$, $\rho_0$& Density and referential density\\
		$g$& Gravity acceleration\\
		$\varphi$ & Generic tracer ($\Theta$ or $S$)\\
		$z$, $z^s$& Vertical coordinate and the z-location of the top sea surface\\
		$\nu_h$, $\nu_v$& Horizontal and vertical viscosities\\
		$\kappa_h$, $\kappa_v$& Tracer diffusion coeffficients\\
		$\vec\bk$&  Unit vector pointing in the local vertical direction \\
		$f$& Coriolis parameter\\
		$f_{\text{eos}}$& Equation of state\\
		\hline 
	\end{tabular}
	\caption{List of variables and other entities for the primitive equations \eqref{layerthickness:single}-\eqref{eqstate:single}.}\label{variable:def}
\end{table}

{\color{black}
	To consider the stratification effects, let us partition the vertical range into $L$ segments (each segment represents one layer) and discretize the above primitive equations in the vertical direction
	as done in \cite{ringler2013multi,calandrini2020exponential}.} Let $\phi_k$ be the vertical average of the generic variable $\phi$ in the layer $k$ for $k=1,2,\cdots,L$. In particular, the horizontal velocity $\bu_k$, {\color{black}the layer thickness $h_k$} and the generic tracer variables $\varphi_k$ are all placed in the middle of the layer while  the vertical velocity $w_k$  in the top of each layer, i.e., in a vertically staggered pattern. We adopt similar but simpler notations 
for  some operators in the $z$-direction \cite[Appendix A.3]{ringler2013multi} as follows:{\color{black}
\begin{equation*}
\left\{\begin{split}
&\overline{\lb \phi\rb}_{k}^{m}=(\phi_{k}+\phi_{k+1})/2,\quad
\overline{\lb \phi\rb}_{k}^{t}=(\phi_{k-1}+\phi_{k})/2,\\
&D_{z^{m}}(\phi)_{k}=\frac{\phi_{k}-\phi_{k+1}}{h_k},\quad
D_{z^{t}}(\phi)_{k}=\frac{\phi_{k-1}-\phi_{k}}{\overline{(h)}_{k}^{t}},
\end{split}
\right.
\end{equation*} 
Let us discretize the vertical derivative terms as
\begin{equation*}
\begin{split}
\frac{\partial}{\partial z}(h_k w_k)\approx  D_{z^m}(w)_{k}h_k,\quad
\ w_k\frac{\partial\bu_k}{\partial z}\approx \overline{\lb w D_{z^t}(\bu)\rb}_{k}^{m},\quad
\frac{\partial}{\partial z}({h_k}\varphi_k w_k)\approx D_{z^m}(\overline{\lb \varphi\rb}^{t}w)_kh_k,
\end{split}
\end{equation*}
and the vertical diffusion terms as
\begin{equation*}
\begin{split}
\frac{\partial}{\partial z}(\nu_v\frac{\partial\bu_k}{\partial z})\approx D_{z^{m}}(\nu_v D_{z^t}(\bu))_k,\quad
\frac{\partial}{\partial z}(\kappa_v\frac{\partial \varphi_k}{\partial z})\approx D_{z_{m}}(\kappa_v D_{z^t}(\varphi))_k.
\end{split}
\end{equation*}}
Then we obtain the multi-layer primitive equations \cite{ringler2013multi}: for $k=1,\dots,L,$
	\begin{empheq}[left={\empheqlbrace}]{alignat=2}
	&\frac{\partial h_k}{\partial t}+\nabla\cdot(h_k\bu_k)
		+D_{z^m}(w)_{k}h_k
	=0,\label{layerthickness}
	\\
	&\frac{\partial \bu_k}{\partial t}+\frac{1}{2}\nabla |\bu_k|^2+(\vec\bk\cdot\nabla\times\bu_k)\bu_k^{\perp}+f\bu_k^{\perp}
	{\color{black}
		+\overline{\lb w D_{z^t}(\bu)\rb}_{k}^{m}
	}
	=-\frac{1}{\rho_0}\nabla p_k+\nu_h\nabla^2\bu_k\nonumber\\
	&\qquad\qquad\qquad
	{\color{black}
		+D_{z^{m}}(\nu_v D_{z^t}(\bu))_k
	}
	,\label{layer:momentum}
	\\
	&\frac{\partial h_k\varphi_k}{\partial t}+\nabla\cdot(h_k\varphi_k\bu_k)
	{\color{black}
		+D_{z^m}(\overline{\lb \varphi\rb}^{t}w)_kh_k
	}
	=\nabla\cdot(h_k\kappa_h\nabla\varphi_k){\color{black}
		+h_kD_{z_{m}}(\kappa_v D_{z^t}(\varphi))_k
	}
	,\label{tracer}
	\\
	&
	p_k=p_k^s(x,y)+
	{\color{black}
		\sum\limits_{l=1}^{k-1}\rho_lgh_l+\frac{1}{2}\rho_kgh_k,\label{hydrocon}
	}
	\\
	&
	\rho_k=f_{\text{eos}}(\Theta_k,S_k,p_k).\label{eqstate}
	\end{empheq}

{\color{black}
The computation of $w_{k}$ via \eqref{layerthickness} is highly dependent on the chosen vertical coordinate  \cite{petersen2015evaluation}. For example, $w_{k}$ is set to zero in the idealized isopycnal vertical coordinate since there is no vertical transport. For the $z$-level vertical coordinate, all layers have a fixed thickness except for the top one ($k=1$). For the $z$-star vertical coordinate, the layer thickness is proportional to the SSH, hence $w_{k}$ is non-zero for all the layers when the ocean is not at rest. 
Petersen et al. \cite{petersen2015evaluation} proposed the arbitrary Lagrangian-Eulerian (ALE) vertical coordinates to obtain $w_{k}$ in a unified way, in which  a new quantity, named $h_{k}^{\text{ALE}}$, is introduced to represent the expected  layer thickness. 
	In this paper, we consider the $z$-star vertical coordinate approach in our method and numerical experiments, which is also the default setting in MPAS-Ocean. 
	}

	\begin{remark}\label{lap-bihar}
		In the momentum equation \eqref{layer:momentum} and \eqref{layer:momentum:single}, a Laplacian form of the horizontal diffusion is present. This is often replaced by alternative formulations or turbulent closures for horizontal dissipation, for example, the hyperviscosity $-\nu_h\nabla^4\bu_k$ will be considered instead of $\nu_h\nabla^2\bu_k$ in the global ocean circulation because it can provide a more scale-selective dissipation effect. 	
		\end{remark}
			
	In the following, we will concentrate on designing high-order time-stepping schemes for the ocean dynamics part of the primitive equations, i.e., the thickness and momentum equations \eqref{layerthickness}-\eqref{layer:momentum} coupled with the hydrostatic condition \eqref{hydrocon}, 
	which together are the key component in numerical ocean modeling. 
	Because the dynamics involves scales having significantly different characteristic times, it is a widespread practice to recognize the fast and slow motions, and split the momentum equation into barotropic and baroclinic subsystems. Define the barotropic and baroclinic velocities by $\overline{\bu}$ and $\widetilde\bu_k$, respectively, by \cite{ringler2013multi,lan2021parallel}:
	\begin{equation}\label{BBspl}
	\left\{
	\begin{split}
	\overline{\bu}&=\sum_{k=1}^{L}h_k\bu_k\Big/\sum_{k=1}^{L}h_k, \\
	\widetilde\bu_k&=\bu_k-\overline{\bu},\qquad k=1,\dots, L.
	\end{split}
	\right.
	\end{equation}
	For simplicity of exposition, we also set the quantities $\bm T_{k}^{u}$, $\bm D_{k}^{u}$ and $\bm T_{k}^{h}$ to be 
	\begin{equation*}
	\begin{split}
	\bm T_{k}^{u}(\bu)&=-\frac{1}{2}\nabla |\bu_k|^2-(\vec\bk\cdot\nabla\times\bu_k)\bu_k^{\perp}-
		\overline{\lb w D_{z^t}(\bu)\rb}_{k}^{m}-\frac{1}{\rho_0}\nabla p_k+\nu_h\nabla^2\bu_k,\\
	\bm D_{k}^{u}(\bu)&=
		D_{z^{m}}(\nu_v D_{z^t}(\bu))_k,\\
	\bm T_{k}^{h}(h, \bu)&=-\nabla\cdot(h_k\bu_k)
		-D_{z^m}(w)_{k}h_k.
	\end{split}
	\end{equation*} 
	Taking the layer-thickness-weighted average of (\ref{layer:momentum}), one then can reformulate the dynamics of the primitive equations into the following baroclinic-barotropic split system:	
	\begin{subequations}\label{bcl:btr:split}
	\begin{align}
	\text{\rm\bf Baroclinic mode:}\qquad&\frac{\partial\widetilde\bu_k}{\partial t}=-f\widetilde\bu_{k}^{\perp}+\bm T_{k}^{u}(\bu)+\bm D_{k}^{u}(\bu)+g\nabla\zeta-\overline{\bm G}, \label{eqclinic}
	\\[4pt]
	\text{\rm\bf Barotropic mode:}\qquad&\frac{\partial}{\partial t}\left(\begin{array}{c}
	\overline\bu\\
	\zeta
	\end{array}\right)
	=-\left(\begin{array}{c}
	f\vec\bk\times \overline\bu+g\nabla\zeta\\
	\nabla\cdot\lb\overline\bu(\zeta+H)\rb
	\end{array}\right)+
	\left(\begin{array}{c}
	\overline{\bm G}\\
	0
	\end{array}\right),\label{eqtropic}
	\\[4pt]
	\text{\rm\bf Fluid thickness:}\qquad&\frac{\partial h_k}{\partial t}=\bm T_{k}^{h}(h, \bu),\label{thickness}
	\end{align}
	\end{subequations}
	where $\zeta$ is the sea surface height (SSH) perturbation, $H$ is the total column height (with respect to the SSH), and the barotropic forcing term $\overline{\bm G}$ contains all other nonlinear terms in the equations for the baroclinic and barotropic modes. It is easy to see that the baroclinic mode subsystem \eqref{eqclinic} and the fluid thickness equations \eqref{thickness} are 
  three-dimensional problems because they involve all $L$ layers whereas the barotropic subsystem \eqref{eqtropic} is only a two-dimensional problem. 
	Because the explicit formula of $\overline{\bm G}$ is hard to present,  an elegant forward-Euler method was proposed in \cite{higdon2002two} for its numerical evaluation  by imposing the layer-thickness-weighted average free condition on the baroclinic mode, which avoids using the explicit expression of  $\overline{\bm G}$. It is worth noting that the SSH perturbation $\zeta$ can be determined from the layer thickness by using the definition $\zeta=\sum_{k=1}^{L}h_k-H$, but it also can be solved from the barotropic subsystem \eqref{eqtropic}. At the space (horizontal) and time continuous  level, these two ways are obviously equivalent, {\color{black} which  implies that the fast time scale is also implicitly 
included in the layer thickness to some extent}.

\section{High-order multirate explicit time-stepping schemes for the split system}\label{sect:ssprk}

In order to  efficiently and accurately simulate the ocean dynamics, we proposed two high-order multirate explicit time-stepping methods for solving the baroclinic-barotropic split system \eqref{bcl:btr:split}. Our method  is based on the framework of the SSPRK method \cite{gottlieb2001strong, gottlieb2011strong} with special treatments of the coupling of the fast and slow modes. In the following, we first review some classic SSPRK schemes.
	Consider the system of ODEs
	\begin{equation}
		 \partial_{t} \bV = F(\bV),
	\end{equation}
where $\bV(t)$ is a vector of unknown variables with $\bV^0=\bV(0)$.
Given  $\bV^n$ at the time $t_{n}$ and the time-step size $\Delta t$, to find $\bV^{n+1}$ at the time $t_{n+1} = t_n+\Delta t$, the second and third-order SSPRK schemes (referred to as SSPRK2 and SSPRK3, respectively) use the following convex combinations of first-order accurate forward-Euler steps:
	\begin{itemize}
		\item SSPRK2 (two stages)
		\begin{equation}\label{SSPRK2}
		\left\{\begin{split}
		\widehat{\bV}^{n+1} &=\bV^{n}+\Delta t F(\bV^{n}), \vspace{0.1cm}\\
		\bV^{n+1}&=\textstyle\frac{1}{2} \bV^{n}+ \frac{1}{2} \left (\widehat{\bV}^{n+1} +\Delta t F(\widehat{\bV}^{n+1})\right). 
		\end{split}\right.
		\end{equation}
		\item SSPRK3 (three stages)
		\begin{equation}\label{SSPRK3}
		\left\{
		\begin{split}
		\widehat{\bV}^{n+1} &=\bV^{n}+\Delta t F(\bV^{n}), \vspace{0.1cm}\\
		\widehat{\bV}^{n+\ahalf}&=\textstyle\frac{3}{4} \bV^{n}+ \frac{1}{4} \left (\widehat{\bV}^{n+1} +\Delta t F(\widehat{\bV}^{n+1})\right ),\vspace{0.1cm}\\
		\bV^{n+1}&=\textstyle\frac{1}{3} \bV^{n}+ \frac{2}{3} \left (\widehat{\bV}^{n+\ahalf} +\Delta t F(\widehat{\bV}^{n+\ahalf})\right ).
		\end{split}
		\right.
		\end{equation}
	\end{itemize}
	We will fit the baroclinic-barotropic split system (\ref{bcl:btr:split}) into the above SSPRK framework with the goal being to develop high-order multirate explicit time-stepping schemes that could stably evolve the baroclinic mode with a large time step and the barotropic mode with a small time step. 

\subsection{Forward-Euler stepping for the baroclinic mode}

	The baroclinic and barotropic modes represent the two different time scales of the ocean dynamics  and thus correspond to different CFL conditions. The former usually admits a large time-step size for stepping, whereas the latter generally requires a small time-step size, but we must take careful care of the coupling between these two modes in order to design  numerical schemes with high-order accuracy. In particular, the barotropic forcing term $\overline{\bm G}$ plays a vital role in their coupling. Without using the  explicit formula of $\overline{\bm G}$, it was  proposed in \cite{higdon2005two} that an effective way to numerically calculate this term  when the forward-Euler stepping is used for solving the baroclinic subsystem \eqref{eqclinic}; this approach is summarized in Algorithm \ref{solveBcl}.  In the sequel, we refer this process as:
	\begin{equation}
	\left[{\widetilde{\bu}}^{1}, \overline{\bm G}\right]=\text{Baroclinic\_FEuler}\lb \bu, \widetilde\bu, \zeta, h,\Delta t\rb.
	\end{equation}  
	
	{\LinesNumberedHidden
		\begin{algorithm}[!ht] \small
			\DontPrintSemicolon
			\KwInput{$\bu$, $\widetilde\bu$, $\zeta$, $h$, $\Delta t$}
			\KwOutput{${\widetilde{\bu}}^{1}$, $\overline{\bm G}$}
			${\widetilde{\bu}}_{k}^{1}=\widetilde\bu_{k}+\Delta t\lb -f \widetilde\bu_{k}^{\perp}+\bm T_{k}^{u}(\bu)+\bm D_{k}^{u}(\bu)+g\nabla\zeta\rb,\ k=1,\dots,L$\;
			$\overline{\bm G}=\frac{1}{\Delta t}\sum_{k=1}^{L}h_{k}{\widetilde{\bu}}_{k}^{1}\Big/\sum_{k=1}^{L}h_{k}$\; 			${\widetilde{\bu}}_{k}^{1}={\widetilde{\bu}}_{k}^{1}-\Delta t\overline{\bm G},\ k=1,\dots,L$
			\caption{ Baroclinic\_FEuler}\label{solveBcl}
		\end{algorithm}
	}
	
				In some ocean modeling situations,  the vertical motion is stiffer than the horizontal motion due to  the much smaller vertical mesh sizes compared with 
		the horizontal ones. A ``vertical mixing'' technique could be applied to tackle this issue, treating the vertical diffusion term $D_{k}^{u}(\bu)$ implicitly in \eqref{eqclinic} (essentially a first-order operator splitting technique). Considering the boundary conditions, such as the bottom drag, we perform the vertical mixing on the whole velocity $\bu$ instead of only on the baroclinic velocity $\widetilde{\bu}$; our approach is summarized in Algorithm \ref{solveBcl_implicit} and is referred  as		
		\begin{equation}
		\left[{\widetilde{\bu}}^{1}, \overline{\bm G}\right]=\text{Baroclinic\_FEuler\_Mixing}\lb \bu, \widetilde\bu, \overline\bu, \zeta, h, \Delta t\rb.
		\end{equation}
		
		{\LinesNumberedHidden
			\begin{algorithm}[!ht] \small
				\DontPrintSemicolon
				\KwInput{$\bu$, $\widetilde\bu$, $\overline\bu$, $\zeta$, $h$, {$\Delta t$}}
				\KwOutput{${\widetilde{\bu}}^{1}$, $\overline{\bm G}$}
				${\widetilde{\bu}}_{k}^{1}=\widetilde\bu_{k}+\Delta t\lb -f \widetilde\bu_{k}^{\perp}+\bm T_{k}^{u}(\bu)+g\nabla\zeta\rb,\ k=1,\dots,L$\;
				$\overline{\bm G}=\frac{1}{\Delta t}\sum_{k=1}^{L}h_{k}{\widetilde{\bu}}_{k}^{1}\Big/\sum_{k=1}^{L}h_{k}$\; 
						${\widetilde{\bu}}_{k}^{1}={\widetilde{\bu}}_{k}^{1}-\Delta t\overline{\bm G},\ k=1,\dots,L$\;
				$\bu^{1}-\Delta t\bm D_{k}^{u}(\bu^{1})={\widetilde{\bu}}_{k}^{1}+\overline\bu$\hspace{2cm}/* {\ttfamily Vertical diffusion solve} */\;
				${\widetilde{\bu}}_{k}^{1}=\bu^{1}-\overline\bu^{n}$
				\caption{ Baroclinic\_FEuler\_Mixing}\label{solveBcl_implicit}
			\end{algorithm}
		}
			
		 When one uses the Forward-Euler scheme combined with vertical mixing for advancing the baroclinic subsystem \eqref{eqclinic} (i.e., Algorithm \ref{solveBcl_implicit} instead of Algorithm \ref{solveBcl}),   the vertical diffusion equation involving $\bu^{1}$ is solved implicitly. The resulting  linear system consists of  blocks of tridiagonal coefficient submatrices  after the spatial discretization, thus it can be efficiently  solved via a tridiagonal matrix algorithm \cite{golub2013matrix} along each independent vertical direction.
 Consequently, the vertical mixing does not appreciably decrease the  performance of explicit time-stepping as shown in the Community Vertical Mixing Project (CVMix) \cite{griffies2015theory}).  On the other hand,  we also note that the
 use of such a treatment could theoretically and practically degrade the high-order accuracies (down to 1 in the worst-case scenario) of the time-stepping schemes  when the effect of the vertical diffusion  term $D_{k}^{u}(\bu)$ on the whole system is very strong. 
 Algorithm \ref{solveBcl} (or Algorithm \ref{solveBcl_implicit} if needed) will be repeatedly used in constructing our SSPRK-based higher-order multirate explicit time-stepping schemes and producing  the needed values of the barotropic forcing term $\overline{\bm G}$ at each stage. 
		
\subsection{SSPRK2-based multirate split-explicit scheme}\label{sec:ssprkse2}

	The proposed SSPRK2-based multirate explicit time-stepping  scheme for solving the baroclinic-barotropic split system \eqref{bcl:btr:split} (referred to as ``SSPRK2-SE'') consists of three steps at each time step from $t_n$ to $t_{n+1}$.
	Two SSPRK2 substepping processes for the barotropic subsystem \eqref{eqtropic} needs to be performed: \textcolor{black}{ the first one is used to predict the intermediate barotropic velocity at $t^{n+1}$}, and the second one is to used to correct and obtain the barotropic velocity. The SSPRK2 substepping is described in Algorithm \ref{sub-step:ssprk2}, where $M$ denotes the number of total substeps (i.e., the small time-step size is set to be $\Delta t/M$). We refer to it as
	\begin{equation}
	\left[ \overline\bu^{1}\right]=\text{Barotropic\_SSPRK2\_Substep}\lb\overline\bu, \zeta, \overline{\bm G}, \Delta t, M\rb.
	\end{equation}
	Let us  define an interpolation operator as
\begin{equation}\label{force2}
		\text{Interp}_2(\overline{\bm G}^1,\overline{\bm G}^2)=
		\frac12\overline{\bm G}^1+\frac12\overline{\bm G}^2.
		\end{equation} 	

	{\LinesNumberedHidden
		\begin{algorithm}[!ht] \small
			\SetKw{KwBy}{by}
			\DontPrintSemicolon
			\KwInput{$\overline\bu$, $\zeta$, $\overline{\bm G}$, $\Delta t$, $M$}
			\KwOutput{$\overline\bu^{1}$} 
			$\overline\bu^{n+0/M}\leftarrow \overline\bu^{n}$, $\zeta^{n+0/M}\leftarrow \zeta$\\ 
			\For{$j= 1$ \KwTo $M$}{
				$\widehat{\overline\bu}^{j/M}=\overline\bu^{(j-1)/M}-\frac{\Delta t}{M} \lb f\vec\bk\times \overline\bu^{(j-1)/M}+g\nabla\zeta^{(j-1)/M}-\overline{\bm G}\rb$\;
				$\widehat{\zeta}^{j/M}=\zeta^{(j-1)/M}-\frac{\Delta t}{M} \nabla\cdot\lb\overline\bu^{(j-1)/M}(\zeta^{(j-1)/M}+H)\rb$\;
				$\widehat{\overline\bu}^{(j+1)/M}=\widehat{\overline\bu}^{j/M}-\frac{\Delta t}{M} \lb f\vec\bk\times \widehat{\overline\bu}^{j/M}+g\nabla\widehat{\zeta}^{j/M}-\overline{\bm G}\rb$\;
				$\widehat{\zeta}^{(j+1)/M}=\widehat{\zeta}^{j/M}-\frac{\Delta t}{M} \nabla\cdot\lb\widehat{\overline\bu}^{j/M}(\widehat{\zeta}^{j/M}+H)\rb$\;
				$\overline\bu^{j/M}=\frac12\lb \overline\bu^{(j-1)/M}+ \widehat{\overline\bu}^{(j+1)/M}\rb$\;
				$\zeta^{j/M}=\frac12\lb \zeta^{(j-1)/M}+ \widehat{\zeta}^{(j+1)/M}\rb$
			}
			\caption{{Barotropic\_SSPRK2\_Substep}}\label{sub-step:ssprk2}
		\end{algorithm}
}
	
 At each time-stepping from $t_n$ to $t_{n+1}$, the proposed  SSPRK2-SE scheme reads as follows: 
	\begin{itemize}
	\item {\bf Preprocessing:} Compute $\overline\bu^n$ and $\widetilde\bu^n$ using the splitting formula \eqref{BBspl}.
	\item {\bf Step 1.} Advance the baroclinic-barotropic system \eqref{bcl:btr:split} from  $\bu^n$ using SSPRK2 with  $\Delta t$ 
	to compute the baroclinic velocity $\widetilde{\bu}^{n+1}$; During the process, the fluid thickness is  also advanced for the first stage in the same way  while one intermediate barotropic velocity is predicted using SSPRK2 subtepping with  $\Delta t/M$ for the purpose of  assistance:
	
	\hspace{0.5cm}{\rm /* {\ttfamily Stage 1 of SSPRK2-SE for baroclinic velocity and fluid thickness} */}
		\begin{subequations}
			\begin{empheq}{align}
			\left[\widehat{\widetilde{\bu}}^{n+1}, \overline{\bm G}^n_0\right]=\text{Baroclinic\_FEuler}\lb \bu^n, \widetilde\bu^n, \zeta^n, h^n,\Delta t\rb.
			\end{empheq}
			{\color{black}\begin{empheq}{align}\label{ssprk2:substep1}
			\left[ \widehat{\overline\bu}^{n+1}\right]&=\text{Barotropic\_SSPRK2\_Substep}\lb\overline\bu^{n}, \zeta^{n}, \overline{\bm G}^n_0, \Delta t,M\rb.
			\end{empheq}}
			\begin{empheq}{align}
			\widehat{\bu}_{k}^{n+1}&= \widehat{\overline\bu}^{n+1}+\widehat{\widetilde{\bu}}_{k}^{n+1}.
			\end{empheq}
			\begin{empheq}[left=\empheqlbrace]{align}\label{ssprk2:thickness1}
			\widehat{h}_{k}^{n+1}&=h_{k}^{n}+\Delta t \bm T_{k}^{h}(h^n, \bu^n),\\[2pt]
			\widehat\zeta^{n+1}&=\textstyle\sum_{k=1}^{L}\widehat{h}_{k}^{n+1}-H.\label{ssprk2zeta1}
			\end{empheq}
		\end{subequations}
		\hspace{0.5cm}{\rm /* {\ttfamily Stage 2 of SSPRK2-SE for baroclinic velocity} */}
		\begin{subequations}
			\begin{empheq}[left=\empheqlbrace]{align}
			\left[\widehat{\widetilde{\bu}}^{n+2}, \overline{\bm G}^{n}_1\right]&=\text{Baroclinic\_FEuler}\lb \widehat\bu^{n+1}, \widehat{\widetilde{\bu}}^{n+1}, \widehat\zeta^{n+1}, \widehat h^{n+1},\Delta t\rb,\\
			\widetilde{\bu}^{n+1}&=\frac12\lb \widetilde{\bu}^n+ \widehat{\widetilde{\bu}}^{n+2}\rb.
			\end{empheq}
			\end{subequations}
			
		\item {\bf Step 2.} Re-advance the barotropic subsystem \eqref{eqtropic} with the interpolated barotropic forcing $\text{Interp}_2(\overline{\bm G}^n_0,\overline{\bm G}^{n}_1)$ from $t_n$ using SSPRK2 substepping with  $\Delta t/M$ to   compute the correct  barotropic velocity $\overline\bu^{n+1}$:\\
		\begin{subequations}
			\begin{empheq}{align}\label{ssprk2:substep2}
			\left[ \overline\bu^{n+1} \right]=\text{Barotropic\_SSPRK2\_Substep}\lb\overline\bu^{n}, \zeta^{n}, \text{Interp}_2(\overline{\bm G}^n_0, \overline{\bm G}^{n}_1), \Delta t,M\rb.
			\end{empheq}
		\begin{empheq}{align}
			\bu_{k}^{n+1}&= \overline{\bu}^{n+1}+\widetilde{\bu}_{k}^{n+1}.
			\end{empheq}	
			\end{subequations}
			
		\item {\bf Step 3.}  Continue to advance the fluid thickness  to obtain $h_{k}^{n+1}$ and update  the SSH perturbation $\zeta^{n+1}$: 
		
		\hspace{0.5cm}{\rm /* {\ttfamily Stage 2 of SSPRK2-SE for fluid thickness} */}
		\begin{subequations}
			\begin{empheq}[left=\empheqlbrace]{align}\label{ssprk2:output:h}
			\widehat{h}_{k}^{n+2}& = \widehat{h}_{k}^{n+1}+\Delta t \bm T_{k}^{h}(\widehat{h}^{n+1},{\bu}^{n+1}),\;\;\;
			h_{k}^{n+1}=\frac12( h_{k}^{n}+\widehat{h}_{k}^{n+2}),\\
			\zeta^{n+1}&=\textstyle\sum_{k=1}^{L}h_{k}^{n+1}-H.\label{ssprk2zeta2}
			\end{empheq}
		\end{subequations}
	\end{itemize}	
	
The total cost of the proposed SSPRK2-SE scheme per time step of size $\Delta t$ mainly consists of 2 forward-Euler baroclinic mode solves (i.e., ${O}(2LN)$ ) and $4M$ forward-Euler barotropic mode solves (i.e., ${O}(4MN)$), where $N$ denotes the number of unknowns per layer. 

\subsection{SSPRK3-based split-explicit scheme}\label{sec:ssprkse3}

	Under the same design principle, the proposed SSPRK3-based multirate explicit time-stepping  scheme for solving the baroclinic-barotropic split system \eqref{bcl:btr:split} (referred to as ``SSPRK3-SE'')  consists of three steps at each time-stepping from $t_n$ to $t_{n+1}$.  
Three SSPRK3 substepping processes for the barotropic subsystem \eqref{eqtropic} needs to be performed: {\color{black} the first two are  to compute the intermediate barotropic velocities 
	at $t^{n+1}$ and $t^{n+2}$}, and the third one is 
 to correct and obtain the final barotropic velocity. We describe such  SSPRK3 substepping  in Algorithm \ref{sub-step:ssprk3}, and refer to it as
	\begin{equation}
	\left[ \overline\bu^{1}\right]=\text{Barotropic\_SSPRK3\_Substep}\lb\overline\bu, \zeta, \overline{\bm G},\Delta t, M\rb.
	\end{equation} 
Let us  also define an interpolation operator as
\begin{equation}\label{force3}
		\text{Interp}_3(\overline{\bm G}^{1},\overline{\bm G}^{2},\overline{\bm G}^{3})=
		\frac{1}{6}\overline{\bm G}^1+\frac{1}{6}\overline{\bm G}^{2}+\frac{2}{3}\overline{\bm G}^{3}.
		\end{equation}
	
	{\LinesNumberedHidden
		\begin{algorithm}[!ht]\small
			\SetKw{KwBy}{by}
			\DontPrintSemicolon
			
			\KwInput{$\overline\bu$, $\zeta$, $\overline{\bm G}$, $\Delta t$, $M$}
			\KwOutput{$\overline\bu^{1}$}
			
			$\overline\bu^{0/M}\leftarrow \overline\bu$, $\zeta^{0/M}\leftarrow \zeta$\\
			\For{$j= 1$ \KwTo $M$}{
				
				$\widehat{\overline\bu}^{j/M}=\overline\bu^{(j-1)/M}-\frac{\Delta t}{M} \lb f\vec\bk\times \overline\bu^{(j-1)/M}+g\nabla\zeta^{(j-1)/M}-\overline{\bm G}\rb$\;
				$\widehat{\zeta}^{j/M}=\zeta^{(j-1)/M}-\frac{\Delta t}{M} \nabla\cdot\lb\overline\bu^{(j-1)/M}(\zeta^{(j-1)/M}+H)\rb$\;
				
				$\widehat{\overline\bu}^{(j+1)/M}=\widehat{\overline\bu}^{j/M}-\frac{\Delta t}{M} \lb f\vec\bk\times \widehat{\overline\bu}^{j/M}+g\nabla\widehat{\zeta}^{j/M}-\overline{\bm G}\rb$\;
				$\widehat{\zeta}^{(j+1)/M}=\widehat{\zeta}^{j/M}-\frac{\Delta t}{M} \nabla\cdot\lb\widehat{\overline\bu}^{j/M}(\widehat{\zeta}^{j/M}+H)\rb$\;
				$\widehat{\overline\bu}^{(j-\sfrac{1}{2})/M}=\frac34\overline\bu^{(j-1)/M}+\frac14\widehat{\overline\bu}^{(j+1)/M}$, 
				$\widehat{\zeta}^{(j-\sfrac{1}{2})/M}=\frac34\zeta^{(j-1)/M}+\frac14\widehat{\zeta}^{(j+1)/M}$\; 
				$\widehat{\overline\bu}^{(j+\sfrac{1}{2})/M}=\overline\bu^{(j-\sfrac{1}{2})/M}-\frac{\Delta t}{M} \lb f\vec\bk\times \overline\bu^{(j-\sfrac{1}{2})/M}+g\nabla\zeta^{(j-\sfrac{1}{2})/M}-\overline{\bm G}\rb$\;
				$\widehat{\zeta}^{(j+\sfrac{1}{2})/M}=\zeta^{(j-\sfrac{1}{2})/M}-\frac{\Delta t}{M} \nabla\cdot\lb\overline\bu^{(j-\sfrac{1}{2})/M}(\zeta^{(j-\sfrac{1}{2})/M}+H)\rb$\;
				$\overline\bu^{j/M}=\frac13\overline\bu^{(j-1)/M}+ \frac23\widehat{\overline\bu}^{(j+\sfrac{1}{2})/M}$, 
				$\zeta^{j/M}=\frac13 \zeta^{(j-1)/M}+ \frac23\widehat{\zeta}^{(j+\sfrac{1}{2})/M}$ 
			}
			\caption{ Barotropic\_SSPRK3\_Substep}\label{sub-step:ssprk3}
		\end{algorithm}
	}

 At each time step from $t_n$ to $t_{n+1}$,  the proposed  SSPRK3-SE scheme reads as follows: 
	\begin{itemize}
	\item {\bf Preprocessing:} Compute $\overline\bu^n$ and $\widetilde\bu^n$ from  $\bu^n$ using the splitting formula \eqref{BBspl}.
	\item {\bf Step 1.} Advance the baroclinic-barotropic system \eqref{bcl:btr:split} from $t_n$ using SSPRK3 with  $\Delta t$ 
	to compute the baroclinic velocity $\widetilde{\bu}^{n+1}$; During the process, the fluid thickness is  also advanced for the first two stages while  two intermediate barotropic velocities is predicted  using SSPRK3 substepping with  $\Delta t/M$ for the purpose of assistance:
	
	\hspace{0.5cm}{\rm /* {\ttfamily Stage 1 of SSPRK3-SE for baroclinic velocity and fluid thickness} */}
		\begin{subequations}
			\begin{empheq}{align}\label{mrigark:part:init}
			\left[\widehat{\widetilde{\bu}}^{n+1}, \overline{\bm G}^n_0\right]=\text{Baroclinic\_FEuler}\lb \bu^n, \widetilde\bu^n, \zeta^n, h^n,\Delta t\rb.
			\end{empheq}
			{\color{black}
			\begin{empheq}{align}\label{ssprk3:substep1}
			\left[ \widehat{\overline\bu}^{n+1}\right]&=\text{Barotropic\_SSPRK3\_Substep}\lb\overline\bu^{n}, \zeta^{n}, \overline{\bm G}^n_0, \Delta t, M\rb.
			\end{empheq}}
			\begin{empheq}{align}
			\widehat{\bu}_{k}^{n+1}&= \widehat{\overline\bu}^{n+1}+\widehat{\widetilde{\bu}}_{k}^{n+1}.
			\end{empheq}
			\begin{empheq}[left=\empheqlbrace]{align}\label{mrigark:part:output:rk3:1}
			\widehat{h}_{k}^{n+1}&=h_{k}^{n}+\Delta t \bm T_{k}^{h}(h^n, \bu^n),\\[2pt]
			\widehat\zeta^{n+1}&=\textstyle\sum_{k=1}^{L}\widehat{h}_{k}^{n+1}-H.\label{ssprk3zeta1}
			\end{empheq}
		\end{subequations}
		\hspace{0.5cm}{\rm /* {\ttfamily Stage 2 of SSPRK3-SE for baroclinic velocity and fluid thickness} */}
		\begin{subequations}
			\begin{empheq}[left=\empheqlbrace]{align}
			\left[\widehat{\widetilde{\bu}}^{n+2}, \overline{\bm G}^{n}_1\right]&=\text{Baroclinic\_FEuler}\lb \widehat\bu^{n+1}, \widehat{\widetilde{\bu}}^{n+1}, \widehat\zeta^{n+1}, \widehat h^{n+1},\Delta t\rb,\\
			\widehat{\widetilde{\bu}}^{n+\sfrac{1}{2}}&=\frac34\widetilde{\bu}^n+\frac14\widehat{\widetilde{\bu}}^{n+2}.
			\end{empheq}
			{\color{black}
			\begin{empheq}[left=\empheqlbrace]{align}
			\widehat{\overline\bu}^{n+2}&=\text{Barotropic\_SSPRK3\_Substep}\lb\widehat{\overline\bu}^{n+1}, \widehat{\zeta}^{n+1}, \overline{\bm G}^{n}_1, \Delta t, M\rb,\label{ssprk3:substep2}\\
			\widehat{\overline\bu}^{n+\sfrac{1}{2}}&=\frac34\overline{\bu}^n +\frac14\widehat{\overline\bu}^{n+2}.\label{convex:stage2:fast}
			\end{empheq}}
			\begin{empheq}{align}
			\widehat\bu_{k}^{n+\sfrac{1}{2}}&= \widehat{\overline{\bu}}^{n+\sfrac{1}{2}}+\widehat{\widetilde{\bu}}_{k}^{n+\sfrac{1}{2}}.
			\end{empheq}
			\begin{empheq}[left=\empheqlbrace]{align}\label{mrigark:part:output:rk3:2}
			\widehat h_{k}^{n+2}&=\widehat{h}_{k}^{n+1}+\Delta t \bm T_{k}^{h}(\widehat{h}^{n+1}, \widehat \bu^{n+1}),\;\;\;
			\widehat h_{k}^{n+\sfrac{1}{2}}=\frac34 h_{k}^{n}+\frac14 \widehat h_{k}^{n+2},\\
			\widehat \zeta^{n+\sfrac{1}{2}}&=\textstyle\sum_{k=1}^{L}\widehat h_{k}^{n+\sfrac{1}{2}}-H.\label{ssprk3zeta2}
			\end{empheq}
		\end{subequations}
		\hspace{0.5cm}{\rm /* {\ttfamily Stage 3 of SSPRK3-SE  for baroclinic velocity} */}
		\begin{subequations}
			\begin{empheq}[left=\empheqlbrace]{align}
			\left[\widehat{\widetilde{\bu}}^{n+\sfrac{3}{2}}, \overline{\bm G}^{n}_{\sfrac{1}{2}}\right]&=\text{Baroclinic\_FEuler}\lb \widehat\bu^{n+\sfrac{1}{2}}, \widehat{\widetilde{\bu}}^{n+\sfrac{1}{2}}, \widehat\zeta^{n+\sfrac{1}{2}}, \widehat h^{n+\sfrac{1}{2}},\Delta t\rb,\\
			\widetilde{\bu}^{n+1}&=\frac13\widetilde{\bu}^n+\frac23\widehat{\widetilde{\bu}}^{n+\sfrac{3}{2}}.
			\end{empheq}
			\end{subequations}
			
			\item {\bf Step 2.} Re-advance the barotropic subsystem \eqref{eqtropic} with the interpolated barotropic force  $\text{Interp}_3(\overline{\bm G}^n_0,\overline{\bm G}^n_1,\overline{\bm G}^n_{\sfrac{1}{2}})$ from $t_n$ using SSPRK3 substepping with  $\Delta t/M$ to compute the correct barotropic velocity $\overline\bu^{n+1}$:
			\begin{subequations}
			\begin{empheq}{align}\label{ssprk3:substep3}
			\left[ \overline\bu^{n+1}\right]=\text{Barotropic\_SSPRK3\_Substep}\lb\overline\bu^{n}, \zeta^{n},\text{Interp}_3(\overline{\bm G}^n_0, \overline{\bm G}^{n}_1, \overline{\bm G}^{n}_{\sfrac{1}{2}}), \Delta t, M\rb.
			\end{empheq}
			\begin{empheq}{align}
			\bu_{k}^{n+1}&= \overline{\bu}^{n+1}+\widetilde{\bu}_{k}^{n+1}.
			\end{empheq}
			\end{subequations}
			
			\item {\bf Step 3.} Continue to advance the fluid thickness (left part from Step 1) to obtain $h_{k}^{n+1}$ and update  the SSH perturbation $\zeta^{n+1}$: 
					
			\hspace{0.5cm}{\rm /* {\ttfamily Stage 3 of SSPRK3-SE for fluid thickness} */}
			\begin{subequations}
			\begin{empheq}[left=\empheqlbrace]{align}\label{ssprk3:output:h}%
			\widehat h_{k}^{n+\sfrac{3}{2}}&=\widehat{h}_{k}^{n+\sfrac{1}{2}}+\Delta t \bm T_{k}^{h}\lb\widehat{h}^{n+\sfrac{1}{2}}, \frac12(\bu^n+\bu^{n+1})\rb,\;\;\;
			h_{k}^{n+1}=\frac13 h_{k}^{n}+\frac23 \widehat h_{k}^{n+\sfrac{3}{2}},\\
			\zeta^{n+1}&=\textstyle\sum_{k=1}^{L}h_{k}^{n+1}-H.\label{ssprk3zeta3}
			\end{empheq}
		\end{subequations}
	\end{itemize}
	The total cost of the proposed SSPRK3-SE scheme per time step of size $\Delta t$ mainly consists of 3  forward-Euler baroclinic mode solves (i.e., ${O}(3LN)$ ) and $9M$ forward-Euler barotropic mode solves (i.e., ${O}(9MN)$). \textcolor{black}{It is approximately 1.5 times (when $L>>M$, i.e.,  the baroclinic mode solves strongly dominate) to 2.25 times  (when $L<<M$, i.e., the barotropic mode solves strongly dominate) of the cost of the SSPRK2-SE scheme.}
	
 \begin{remark}\label{InconsisError}	
Note that there are two approaches to determine the  SSH perturbation $\zeta$  due the mode splitting, however their equivalence  will not hold anymore in the discrete setting  which could make $\zeta$ become over-determined.  The  inconsistency caused by the determination of the SSH perturbation $\zeta$  after discretizing the 
baroclinic-barotropic split system in time and space  is often referred to as the mode-splitting error in the literature.
\textcolor{black}{In the proposed SSPRK-SE schemes,
we update $\zeta$  using the sum of all layer fluid thickness $\{h_k\}$  (i.e., $\zeta=\sum_{k=1}^{L}h_{k}-H$), which could be different from that produced by the corresponding barotropic substeppings for \eqref{eqtropic} (the fast mode solve). This could consequently result in  additional numerical inconsistency and stability issue since 
 $\{h_k\}$ implicitly contain the fast time scale and need to be consistent with $\zeta$ produced from the fast mode solve}. In order to relieve 
this issue we will incorporate the velocity adjustments for updating the layer fluid thickness (see Subsection \ref{reconciliation}) to reconcile the SSH, which makes the resulting $\zeta$  be the same as that from the corresponding fast mode solve.
	\end{remark}

\subsection{Sea surface height reconciliation in flux-form}\label{reconciliation} 

	As aforementioned, there are two ways to compute the SSH perturbation $\zeta$ due the mode splitting: one is from the barotropic substepping and the other from the layer thickness. To resolve the resulting model inconsistency errors  and further stabilize the fast mode implicitly included in the layer thickness. We propose  is to use a SSH reconciliation process in flux-form  \cite{higdon2005two} within the framework of the proposed SSPRK-SE schemes, i.e., make the value of $\zeta$ computed from the layer fluid  thickness  match that produced from the barotropic substepping at each stage.

\paragraph{\bf SSH reconciliation for SSPRK2-SE}
			
		 For the SSPRK2-SE scheme, \textcolor{black}{the SSH perturbations obtained from each of the two barotropic SSPRK2 substeppings (\eqref{ssprk2:substep1} and \eqref{ssprk2:substep2}, respectively) can be expressed in the following representation: for $i=1,2$,}
		\begin{equation*}
		\zeta_i=\zeta^{n}-\Delta t \nabla\cdot \bF_i^{\zeta}
		\end{equation*}
		with the accumulated flux
		 $$\bF_i^{\zeta}=\sum_{j=1}^{M}\frac{\overline\bu^{n+(j-1)/M}(\zeta^{n+(j-1)/M}+H)+\widehat{\overline\bu}^{n+j/M}(\widehat{\zeta}^{n+j/M}+H)}{2M},$$ where $i=1$ is for  \eqref{ssprk2:substep1}  and $i=2$ for \eqref{ssprk2:substep2}. 
		{Here for simplicity we use the same notations for the intermediate substep velocities and SSH perturbations at both cases, but their values may differ.} On the other hand, the corresponding SSH perturbations (\eqref{ssprk2zeta1} and \eqref{ssprk2zeta2}) obtained by summing up the layer-thickness are given by: for $i=1,2$, 
		\begin{equation*}
		\zeta_i= \zeta^{n}-\Delta t \nabla\cdot \bF_{i}^{h},
		\end{equation*}
		with{\color{black}
		\begin{equation*}
		\bF_{i}^{h}=\left\{\begin{array}{ll}
		\sum_{k=1}^{L} \bu_{k}^{n}h_{k}^{n},&i=1,\\[6pt]
		\sum_{k=1}^{L} \frac{\bu_{k}^{n}h_{k}^{n} + \bu_{k}^{n+1}\widehat{h}_{k}^{n+1}}{2},&i=2.  		\end{array}
		\right.
		\end{equation*}}
		Thus the corresponding flux deficits between each pair of them are given by 	
		$D_i^e=\bF_{i}^{\zeta}-\bF_{i}^{h},\ i=1,2$, and  we need to respectively compensate them back by slightly modifying 
		{the calculations of the layer thickness at each stage}.
For the efficiency of the computation,
		we will adjust the layer thickness flux by explicitly assigning the deficit $D_i^e$ back in a layer thickness-weighted way. Hence, we define the transport velocity adjustments{\color{black}
	\begin{equation}
		\bu_i^{A,n+1}=\left\{
		\begin{array}{ll}
				\frac{D_1^e}{\sum_{k=1}^{L} h_{k}^{n}} ,&i=1,\\[10pt]
				\frac{2D_2^e}{\sum_{k=1}^{L} \widehat{h}_{k}^{n+1}} ,&i=2,
		\end{array}
		\right.
		\end{equation}}
		where the factor 2 when $i=2$ is due  to the final averaging step of layer thickness in  SSPRK2-SE. Finally, our SSH reconciliation process  is completed by revising the layer thickness updates \eqref{ssprk2:thickness1} and \eqref{ssprk2:output:h} in SSPRK2-SE  to
		{\color{black}\begin{equation}
		\widehat{h}_{k}^{n+1}=h_{k}^{n}+\Delta t \bm T_{k}^{h}(h^n, \bu^n+\bu_1^{A,n+1})
		\end{equation}}
		and 
		\begin{equation}
		\widehat{h}_{k}^{n+2} = \widehat{h}_{k}^{n+1}+\Delta t \bm T_{k}^{h}(\widehat{h}^{n+1},{\bu}^{n+1}+\bu_2^{A,n+1}),\;\;\;h_{k}^{n+1}=\frac12( h_{k}^{n}+\widehat{h}_{k}^{n+2})
		\end{equation}
		respectively.
		
		\paragraph{\bf SSH reconciliation for SSPRK3-SE}
		For the SSPRK3-SE scheme, \textcolor{black}{the SSH perturbation obtained from each of the three barotropic SSPRK3 substeppings (\eqref{ssprk3:substep1}, \eqref{ssprk3:substep2} and \eqref{ssprk3:substep3}, respectively) is: {for $i=1, 2, \text{ and } 3,$}}
		\begin{equation*}
		\zeta_i=\zeta^{\delta_i}-\Delta t \nabla\cdot \bF_{i}^{\zeta},
		\end{equation*}
		where \textcolor{black}{$\zeta^{\delta_1}=\zeta^{\delta_3}=\zeta^{n}$, $\zeta^{\delta_2}=\widehat{\zeta}^{n+1}$,} and 
		$$\bF_{i}^{\zeta}=\sum_{j=1}^{M}\lb \frac{\overline\bu^{(j-1)/M}(\zeta^{(j-1)/M}+H)+\widehat{\overline\bu}^{j/M}(\widehat{\zeta}^{j/M}+H)}{6M}+\frac{2\overline\bu^{(j-\sfrac{1}{2})/M}(\zeta^{(j-\sfrac{1}{2})/M}+H)}{3M} \rb$$
is the accumulated flux. Note that for $i=2$,  we need reconcile SSH perturbation at $t_{n+2}$ instead of $t_{n+\frac12}$ since we only have its value at $t_{n+2}$ after the substepping in \eqref{ssprk3:substep2} and do not further take the convex combination in \eqref{convex:stage2:fast}. Similarly, we use the same notations for the intermediate substep velocities and SSH perturbations for all of them but their values may differ at the different places. 		
		 On the other hand, the corresponding SSH perturbations obtained by summing up the layer thicknesses are given by: for $i=1, 2, \text{ and } 3,$
		\begin{equation*}
\zeta_i=\zeta^{\delta_i}-\Delta t \nabla\cdot \bF_{i}^{h},
\end{equation*}
			with{\color{black}
		\begin{equation*}
		\bF_{i}^{h}=\left\{
		\begin{array}{ll}
		\sum_{k=1}^{L}\bu_k^nh_k^n,&i=1,\\[6pt]
		\sum_{k=1}^{L}\widehat \bu_{k}^{n+1}\widehat{h}_{k}^{n+1},&i=2,\\[6pt]
		\sum_{k=1}^{L}\frac{\bu_k^nh_k^n+\widehat \bu_{k}^{n+1}\widehat{h}_{k}^{n+1}+ 2(\bu_{k}^n+\bu_{k}^{n+1})\widehat{h}_{k}^{n+\sfrac{1}{2}}}{6},&i=3.
		\end{array}
		\right.
		\end{equation*}}
		Thus the corresponding flux deficits are given by 	
		$D_i^e=\bF_{i}^{\zeta}-\bF_{i}^{h},\ i=1, 2, \text{ and }3$.
		Following the same process for SSPRK2-SE, we will adjust the layer thickness flux by explicitly assigning the deficit  $D_i^e$ back in a layer thickness-weighted way. The transport velocity adjustments are defined by {\color{black}
		\begin{equation}
		\bu_i^{A,n+1}=\left\{
		\begin{array}{ll}
		\frac{D_1^e}{\sum_{k=1}^{L} h_k^n},&i=1,\\[10pt]
		\frac{D_2^e}{\sum_{k=1}^{L} \widehat{h}_{k}^{n+1}},&i=2,\\[10pt]
		\frac{\frac32D_3^e}{\sum_{k=1}^{L} \widehat{h}_{k}^{n+\sfrac{1}{2}}} ,&i=3,
		\end{array}
		\right.
		\end{equation}}
		where the factor $\frac{3}{2}$ when $i=3$ is due  to the final averaging step of layer thickness in  SSPRK3-SE. Finally, our SSH reconciliation process  is completed by revising the layer thickness \eqref{mrigark:part:output:rk3:1}, \eqref{mrigark:part:output:rk3:2} and \eqref{ssprk3:output:h} in SSPRK3-SE to
		{\color{black}\begin{equation}
		\widehat{h}_{k}^{n+1}=h_{k}^{n}+\Delta t \bm T_{k}^{h}(h^n, \bu^n+\bu_1^{A,n+1}),
		\end{equation}}
		{\color{black}\begin{equation}
		\widehat h_{k}^{n+2}=\widehat{h}_{k}^{n+1}+\Delta t \bm T_{k}^{h}(\widehat{h}^{n+1}, \widehat \bu^{n+1}+\bu_2^{A,n+1}),\;\;\;
		\widehat h_{k}^{n+\sfrac{1}{2}}=\frac34 h_{k}^{n}+\frac14 \widehat h_{k}^{n+2},
		\end{equation}}
		and
		\begin{equation}
		\widehat h_{k}^{n+3/2}=\widehat{h}_{k}^{n+\sfrac{1}{2}}+\Delta t \bm T_{k}^{h}\lb\widehat{h}^{n+\sfrac{1}{2}}, \frac12(\bu_{k}^n+\bu_{k}^{n+1})+\bu_3^{A,n+1}\rb,\;\;\;
		h_{k}^{n+1}=\frac13 h_{k}^{n}+\frac23 \widehat h_{k}^{n+\sfrac{3}{2}}.
		\end{equation}
		respectively.
	
	\begin{remark}
	With the SSH reconciliation processes to further remove the model inconsistency errors due to the mode splitting, it is expected the  numerical stability of the proposed SSPRK2-SE 
	scheme follows that of SSPRK2 with $\Delta t$ for the baroclinic mode solve and   $\Delta t /M$ for the barotropic mode solve, and 
	so does the SSPRK3-SE scheme in correspondence with SSPRK3.  \textcolor{black}{In addition, since that SSPRK3 allows even larger  time-step sizes (i.e.,  CFL conditions) compared to SSPRK2 \cite{kubatko2014optimal,kubatko2008time},  we also expect that 
	the SSPRK3-SE scheme is  numerically more stable than the SSPRK2-SE scheme, which will be verified through experiments in Section 	\ref{sect:experiments}.}
\end{remark}	
	
	\section{Temporal error analysis}\label{sec:err}

	In this section, we  analyze and discuss the temporal errors  of the proposed SSPRK-SE schemes.
	Note that in the following analysis we only consider the semi discrete-in-time case of the baroclinic-barotropic split system \eqref{bcl:btr:split}, i.e., spatial discretization is not considered.
	It is also assumed that neither vertical mixing nor SSH reconciliation is  applied  and the solutions possess sufficient smoothness. To simplify the discussion, we first define $\overline{\bm G}^*$ as the general barotropic forcing term
	 and introduce the following extra notations:
	\begin{align*}
	&F^n=f\vec\bk\times \overline\bu^{n}+g\nabla\zeta^{n}-\overline{\bm G}^*,\quad
	P^n=\nabla\cdot\lb \overline\bu^{n}(\zeta^{n}+H)\rb,\\
	&E^n=f\vec\bk\times F^n+g\nabla P^n,\quad
	Q^n=\nabla\cdot\lb \bu^{n}P^n+(\zeta^{n}+H)F^n\rb.
	\end{align*}
Let us first consider Algorithms \ref{sub-step:ssprk2} with regards to its  temporal accuracies. 	
	\begin{proposition}\label{prop:ssprk2}
		Given 
		\begin{equation*}
		\left[ \overline\bu^{n+1}\right]=\text{\rm Barotropic\_SSPRK2\_Substep}\lb\overline\bu^{n}, \zeta^{n}, \overline{\bm G}^*, \Delta t, M\rb,
		\end{equation*}
		then when $\Delta t$ sufficiently small,
		\begin{equation}\label{ssprk2:formula}
		\overline\bu^{n+1}=\;\overline\bu^{n}-\Delta tF^n+\frac{\Delta t^2}{2}E^n-\frac{M^2-1}{3M^2}\Delta t^3\lb f\vec\bk\times E^n+g\nabla Q^n\rb+O(\Delta t^4).
		\end{equation}
	\end{proposition}
	\begin{proof}
		Given $\lb \overline\bu^{n+(j-1)/M}, \zeta^{n+(j-1)/M}\rb$,   after one substep in Algorithm \ref{sub-step:ssprk2}, we have
		\begin{equation*}
		\left\{\begin{split}
		\overline\bu^{n+j/M}=\;&\overline\bu^{n+(j-1)/M}-\frac{\Delta t}{M}F^{n+(j-1)/M}+\frac{\Delta t^2}{2M^2}E^{n+(j-1)/M},\\
		\zeta^{n+j/M}=\;&\zeta^{n+(j-1)/M}-\frac{\Delta t}{M}P^{n+(j-1)/M}+\frac{\Delta t^2}{2M^2}Q^{n+(j-1)/M}-\frac{\Delta t^3}{2M^3}\nabla\cdot\lb F^{n+(j-1)/M}P^{n+(j-1)/M}\rb.
		\end{split}
		\right.
		\end{equation*}
		By mathematical induction, we  obtain 
		\begin{equation}\label{ssprk2:formula:two}
		\left\{\begin{split}
		\overline\bu^{n+j/M}&=\overline\bu^{n}-\frac{j\Delta t}{M}F^n+\frac{(j\Delta t)^2}{2M^2}E^n-\frac{j(j^2-1)}{3M^3}\Delta t^3\lb f\vec\bk\times E^n+g\nabla Q^n\rb+O(\Delta t^4),\\
		\zeta^{n+j/M}&=\zeta^{n}-\frac{j\Delta t}{M}P^n+\frac{(j\Delta t)^2}{2M^2}Q^n+O(\Delta t^3).
		\end{split}
		\right.
		\end{equation}
		Then (\ref{ssprk2:formula}) follows directly by setting $j=M$ in \eqref{ssprk2:formula:two}. 
	\end{proof}
	
\textcolor{black}{It is worth noting that the barotropic velocity   $\overline\bu$ is used only as a part of the flux term (through $\bm T_{k}^{u}(\bu)$ and $\bm D_{k}^{u}(\bu)$) for computing  the baroclinic velocity $ \widetilde\bu$ in the forward-Euler stepping (Algorithm \ref{solveBcl}}). Then we can obtain the following error estimate for the 
SSPRK2-SE scheme.
	
	\begin{thm}\label{ssprk2:accuracy}
		The SSPRK2-SE scheme is of the second-order accuracy in time for solving the baroclinic-barotropic split system \eqref{bcl:btr:split}.
	\end{thm}
	
	\begin{proof}
		It is easy to show that if we add $O(\Delta t^3)$ terms to the forward-Euler steps in the classic SSPRK2 scheme, the resulting scheme is still second-order accurate in time, that is, if we have the following scheme
		\begin{equation*}
		\left\{\begin{split}
		\widehat{\bV}^{n+1} &=\bV^{n}+\Delta t F(\bV^{n})+O(\Delta t^3), \vspace{0.1cm}\\
		\bV^{n+1}&=\textstyle\frac{1}{2} \bV^{n}+ \frac{1}{2} \left (\widehat{\bV}^{n+1} +\Delta t F(\widehat{\bV}^{n+1})+O(\Delta t^3)\right),
		\end{split}\right.
		\end{equation*}
		then $\bV^{n+1}$ is of the second-order accuracy. Our proof will be based on  this result. In the sequel, we mark the values obtained by following the classic  SSPRK2 scheme \eqref{SSPRK2} with the subscript ``RK2".
		
		At the first stage, the barotropic mode solve \eqref{ssprk2:substep1} is different from SSPRK2. According to Proposition \ref{prop:ssprk2} with $\overline{\bm G}^*=\overline{\bm G}^n$, we have
		\begin{equation*}
		\widehat{\overline\bu}^{n+1}=~\overline\bu^{n}-\Delta tF^n+O(\Delta t^2)=~ \widehat{\overline\bu}_{\text{RK2}}^{n+1}+O(\Delta t^2),
		\end{equation*}
		which is equivalent to the forward-Euler with the extra term $O(\Delta t^2)$. 
		However, the barotropic velocity affects the baroclinic velocity and layer thickness at Stage 2 as a part of the  flux term (i.e., after multiplying with $\Delta t$. For the baroclinic mode, we then have
		\begin{equation*}
		\begin{split}
	\widehat{\widetilde{\bu}}^{n+2}&=\text{Baroclinic\_FEuler}\lb \widehat\bu_{\text{RK2}}^{n+1}+O(\Delta t^2), \widehat{\widetilde{\bu}}_{\text{RK2}}^{n+1}, \widehat\zeta_{\text{RK2}}^{n+1}, \widehat h_{\text{RK2}}^{n+1},\Delta t\rb=~\widehat{\widetilde{\bu}}_{\text{RK2}}^{n+2}+O(\Delta t^3),\\
	\widetilde{\bu}^{n+1}&=\frac12 ( \widetilde{\bu}^{n}+ \widehat{\widetilde{\bu}}^{n+2}).
	\end{split}
		\end{equation*}
		Therefore, the baroclinic velocity $\widetilde{\bu}^{n+1}$ obtained by SSPRK2-SE is of second-order accuracy. Next, let us turn to the barotropic mode. SSPRK2-SE recomputes the fast mode with the values at $t_n$ and the predicted term $\overline{\bm G}^{n+1}$. Since $\widehat{\widetilde{\bu}}^{n+2}=~\widehat{\widetilde{\bu}}_{\text{RK2}}^{n+2}+O(\Delta t^3)$, it is easy to see that $\overline{\bm G}^{n+1} = ~\overline{\bm G}_{\text{RK2}}^{n+1}+O(\Delta t^3)$.
		
		By taking the time derivative on \eqref{eqtropic}, we get
		\begin{equation}
		\frac{\partial^2 \overline\bu}{\partial t^2}=f\vec\bk\times\lb f\vec\bk\times\overline\bu+g\nabla\zeta-\overline{\bm G}\rb+g\nabla\lb\nabla\cdot(\overline\bu(\zeta+H))\rb+\frac{\partial \overline{\bm G}}{\partial t}.
		\end{equation}
		According to Proposition \ref{prop:ssprk2}, the barotropic velocity at Stage 2 is given by
		\begin{equation}
		\begin{split}
		\overline\bu^{n+1}=\;&\overline\bu^{n}-\Delta tF^n+\frac{\Delta t^2}{2}E^n+o(\Delta t^2)\\
		=\;&\overline\bu^{n}-\Delta t\lb f\vec\bk\times \overline\bu^{n}+g\nabla\zeta^{n}-\frac{1}{2}\lb \overline{\bm G}^n+\overline{\bm G}^{n+1}\rb\rb\\
		&+\frac{\Delta t^2}{2}\lb f\vec\bk\times \lb f\vec\bk\times \overline\bu^{n}+g\nabla\zeta^{n}-\frac{1}{2}\lb \overline{\bm G}^n+\overline{\bm G}^{n+1}\rb\rb \rb\\
		&+\frac{\Delta t^2}{2}\lb g\nabla(\nabla\cdot\lb \overline\bu^{n}(\zeta^{n}+H)\rb) \rb+O(\Delta t^3).
		\end{split}
		\end{equation}
		To identify the truncation error of the SSPRK2-SE scheme, let us replace the numerical solutions by their exact  counterparts, which yields
		\begin{equation}
		\begin{split}
		\overline\bu(t_{n+1})=\;&\overline\bu(t_{n})-\Delta t\lb f\vec\bk\times \overline\bu(t_{n})+g\nabla\zeta(t_{n})-\frac{1}{2}\lb \overline{\bm G}(t_{n})+\overline{\bm G}(t_{n+1})\rb\rb\\
		&+\frac{\Delta t^2}{2}\lb f\vec\bk\times \lb f\vec\bk\times \overline\bu(t_{n})+g\nabla\zeta(t_{n})-\frac{1}{2}\lb \overline{\bm G}(t_{n})+\overline{\bm G}(t_{n+1})\rb\rb \rb\\
		&+\frac{\Delta t^2}{2}\lb g\nabla(\nabla\cdot\lb \overline\bu(t_{n})(\zeta(t_{n})+H)\rb) \rb+O(\Delta t^3)\\
		=\;&\overline\bu(t_{n})-\Delta t\lb f\vec\bk\times \overline\bu(t_{n})+g\nabla\zeta(t_{n})-\overline{\bm G}(t_{n})-\frac{\Delta t}{2}\frac{\partial\overline{\bm G}(t_{n})}{\partial t}\rb\\
		&+\frac{\Delta t^2}{2}\lb f\vec\bk\times \lb f\vec\bk\times \overline\bu(t_{n})+g\nabla\zeta(t_{n})  -\overline{\bm G}(t_{n}) \rb \rb\\
		&+\frac{\Delta t^2}{2}\lb g\nabla(\nabla\cdot\lb \overline\bu(t_{n})(\zeta(t_{n})+H)\rb) \rb+O(\Delta t^3)\\
		=\;&\overline\bu(t_{n})-\Delta t\lb f\vec\bk\times \overline\bu(t_{n})+g\nabla\zeta(t_{n})-\overline{\bm G}(t_{n})\rb+\frac{\Delta t^2}{2}\frac{\partial^2 \overline\bu(t_{n})}{\partial t^2}+O(\Delta t^3).
		\end{split}
		\end{equation}
		Thus the barotropic velocity has a third-order truncation error
		and thus is of second-order accuracy in time.
	 For the layer thickness $h$, let us rewrite its updating in a compact way as
		\begin{equation}
		h^{n+1}=h^n-\frac{\Delta t}{2}\lb \bm T_{k}^{h}(h^n, \bu^n)+ \bm T_{k}^{h}(\widehat{h}^{n+1},{\bu}^{n+1})\rb.
		\end{equation} 
		Notice that here we have used the second-order accurate solution ${\bu}^{n+1}$ instead of the first order accurate velocity $\widehat{\bu}_{\text{RK2}}^{n+1}$. In addition, the difference between ${\bu}^{n+1}$ and $\widehat{\bu}_{\text{RK2}}^{n+1}$ is $O(\Delta t^2)$. Based on the classic result of  SSPRK2, we know that
				\begin{equation}
		\widetilde h_{\text{RK2}}^{n+1}=h^n-\frac{\Delta t}{2}\lb\bm T_{k}^{h}(h^n, \bu^n)+ \bm T_{k}^{h}(\widehat{h}_{\text{RK2}}^{n+1},\widehat{\bu}_{\text{RK2}}^{n+1})\rb,
		\end{equation}
		which gives us a third-order truncation error in time, i.e., $|\widetilde h_{\text{RK2}}^{n+1}-h(t_{n+1})| = O(\Delta t^3)$. Furthermore,
		\begin{equation}
		|h^{n+1}-\widetilde h_{\text{RK2}}^{n+1}|=\Big|\frac{\Delta t}{2}\lb \bm T_{k}^{h}(\widehat{h}^{n+1},{\bu}^{n+1})-\bm T_{k}^{h}(\widehat{h}_{\text{RK2}}^{n+1},\widehat{\bu}_{\text{RK2}}^{n+1})\rb\Big|=O(\Delta t^3).
		\end{equation}
		Noticing that $\widehat{h}^{n+1}=~\widehat{h}_{\text{RK2}}^{n+1}$, thus it holds that $h^{n+1}$ is also of
		 	second-order accuracy in time.
	\end{proof}
		\textcolor{black}{
		As for the SSPRK3-SE scheme,  since the substepping for the barotropic mode is applied at each of the first two stages with a first-order  approximation  of the barotropic forcing  $\overline{\bm G}$ obtained from the corresponding foward-Euler solve for the baroclinic mode,  the overall truncation errors could not reach $O(\Delta t^4)$, instead they are still only $O(\Delta t^3)$ by following a similar analysis as above. Thus the convergence of the SSPRK3-SE scheme could downgrade to  the second order in time, which will be also checked through numerical experiments in Section \ref{sect:experiments}}.
		
\section{Numerical experiments}\label{sect:experiments}

We implement the proposed two SSPRK-based multirate explicit  time-stepping schemes (SSPRK2-SE and SSPRK3-SE) within the framework of   
MPAS-Ocean (using its built-in subroutines). In MPAS-Ocean, the baroclinic-barotropic split system (\ref{bcl:btr:split}) is discretized in space by the TRiSK scheme (a specially-designed finite volume approximation) on unstructured, locally orthogonal dual meshes {\rm \cite{ringler2010unified, thuburn2012framework, thuburn2009numerical}} to ensure many physical properties of the space-continuous system in the spatially discrete setting, such as the conservations of mass, total energy, and potential vorticity.
To  investigate  their accuracy and performance, we use two benchmark test cases  from the MPAS-Ocean platform (Version 7.0) \cite{petersen2018mpas}. 
	In these two cases, the bottom drag is considered  as a bottom boundary condition:
	\begin{equation}
	\lim\limits_{z\rightarrow z_{\text{bot}}}\nu_v \frac{\partial \bu}{\partial z}=c_{\text{drag}}|\bu|\bu,
	\end{equation}
	where $|\bu|$ denotes the magnitude of the velocity, $c_{\text{drag}}$ is the bottom drag coefficient, and $z_{\text{bot}}$ is the z-location of the ocean bottom. In the following numerical experiments, the temperature means the potential temperature. In addition, if the results are obtained  with the vertical mixing, we will call Algorithm \ref{solveBcl_implicit} (``Baroclinic\_FEuler\_Mixing") instead of Algorithm \ref{solveBcl} (``Baroclinic\_FEuler") to advance the baroclinic mode in the proposed SSPRK2-SE and SSPRK3-SE schemes. The time-step size used for  the barotropic substepping (Algorithm \ref{sub-step:ssprk2} or \ref{sub-step:ssprk3})   is denoted by $\Delta_{\text{btr}} = \Delta t/M$. In order to measure the accuracy,  the relative $l_{2}$ error is considered for the velocity $\bu_1$ (i.e. the surface velocity) and the layer thickness $h_1$ (i.e.,  fluid thickness of the top layer) at the terminal time $T$: 
	\begin{equation*}
	\frac{\|\bu_{1}^{t}-\bu_{1}^{r}\|_{2}}{\|\bu_{1}^{r}\|_{2}} \text{\qquad and \qquad}\frac{\|h_{1}^{t}-h_{1}^{r}\|_{2}}{\|h_{1}^{r}\|_{2}},
	\end{equation*}
	where $\bu_{1}^{r}$ and $h_{1}^{r}$ are the reference values, $\bu_{1}^{t}$ and $h_{1}^{t}$ are the testing values, and $\|\cdot\|_{2}$ is the vector $l_{2}$-norm. All numerical experiments are performed on the cluster  ``Cori'' at the National Energy Research Scientific Computing Center (NERSC). In particular, we run our codes on  the ``Haswell" processor nodes for which each node has two 16-core Intel Xeon ``Haswell" processors and 128 GB memory.
		
\subsection{The baroclinic eddies test case with twenty layers}\label{barochan}

	We first consider an ideal test case, namely the baroclinic eddies test case \cite{ringler2013multi,petersen2018mpas}, for the primitive equations with 20 vertical layers  provided by the  MPAS-Ocean platform imported from \cite{ilicak2012spurious}. The domain consists of a horizontally periodic channel of latitudinal extent 440~km and
	longitudinal extent 160~km, with a flat bottom of 1 km vertical depth. The channel is on a f-plane \cite{cushman2011introduction} with the Coriolis parameter $f=1.2\times10^{-4}$ s$^{-1}$, where s denotes seconds. The initial temperature
	decreases downward in the meridional direction. A cosine shaped temperature
	perturbation with a wavelength of 120~km in the zonal direction is used to instigate the baroclinic
	instability. 
	 The horizontal domain is partitioned by a 10-km-resolution SCVT mesh \cite{ju2003scvt} such that each layer contains 3,920 cells, 11,840 edges, and 7,920 vertices.  The horizontal  viscosity is given as $\nu_h=10$ and the vertical one as $\nu_v=\text{1.0E-4}$. For this test case, the effect of vertical mixing is not necessary since the vertical diffusion is almost neglectable.
	 We compute the relative $l_{2}$ errors of the velocity and layer thickness approximations at the terminal time $T=4096$~s with the reference solution generated using the SSPRK3-SE scheme with
	$\Delta t=\Delta_{\text{btr}} t=0.25$~s.
	
	We first test the temporal convergence order  of the current dynamical core (the MPAS-SE scheme with default setting) in MPAS-Ocean.  We also comment out the part of the tracers in the code because in the convergence test we are only interested in the ocean dynamics.  The relative $l_{2}$ errors in the velocity and layer thickness approximations are reported in Table \ref{baroclinic:mpas}, which shows that the MPAS-SE scheme blows up when  $\Delta t = 128$~s and the temporal convergence rates are less than one  for both the velocity and the layer thickness. 
	
	\begin{table}[!ht]
		\centering\footnotesize
		\begin{tabular}{|c|rr|rr|}
			\hline
			\multirow{2}{*}{$\Delta t$-$\Delta_{\text{btr}} t$}& \multicolumn{2}{c|}{Velocity}& \multicolumn{2}{c|}{Layer thickness}\\\cmidrule(lr){2-3}\cmidrule(lr){4-5}
			&Error &Rate &Error &Rate \\ \hline\hline
128-128&N/A&-&N/A&-\\
64-64 &2.422E-02& - &1.500E-05& -\\
32-32 &2.142E-02&0.18&1.332E-05&0.17\\
16-16 &1.721E-02&0.32&1.031E-05&0.37\\
08-08 &1.281E-02&0.43&7.264E-06&0.51\\
04-04 &9.030E-03&0.50&4.833E-06&0.59\\
02-02 &6.008E-03&0.59&3.069E-06&0.66\\
01-01 &3.735E-03&0.69&1.851E-06&0.73\\
			\hline
		\end{tabular}
		\caption{Relative $l_{2}$ errors and convergence rates of $\bu_1$ and $h_1$  for the baroclinic  eddies test case produced by the MPAS-SE scheme.}\label{baroclinic:mpas}
	\end{table}
	
	{\color{black}
	Next, we test the proposed  two SSPRK-SE  schemes for the ocean dynamics in terms of accuracy by fixing $\Delta_{\text{btr}} t =$ 64, 32, 16, 8, 4 and 2~s and taking 
	$M=$ 16, 8, 4, 2 and 1 (i.e., $\Delta t= M\Delta_{\text{btr}} t$).
	We first turn off the vertical mixing, and the results produced by SSPRK2-SE and SSPRK3-SE with or without SSH reconciliation are reported in Tables \ref{baroclinic:m:1} and  \ref{baroclinic:m:1:imp}, respectively. Although the vertical diffusion only has negligible effect in this test case, we still apply the vertical mixing in both proposed SSPRK-SE schemes and rerun all the tests.  It is found as expected that  the vertical mixing does not affect the performance of both  SSPRK-SE schemes for this baroclinic eddies test case, i.e., the produced numerical solutions and errors by SSPRK2-SE and SSPRK3-SE with vertical mixing are almost identical to those reported in Tables \ref{baroclinic:m:1} and \ref{baroclinic:m:1:imp}, and thus not shown here.

	For SSPRK2-SE without SSH reconciliation, as $\Delta_{\text{btr}}t$  are uniformly refined from 64~s to 2~s, we observe from Table \ref{baroclinic:m:1} that the optimal second-order temporal  convergence  is gradually achieved for both velocity and fluid thickness  when $M$ is small (such as $M=$ 1, 2 and 4). Meanwhile, the larger M it uses, the quicker the scheme becomes unstable and blows up along the increase of $\Delta_{\text{btr}}t$. On the other hand, it is easy to find that SSPRK2-SE with SSH reconciliation is  more stable and accurate when $M$ is large, and the optimal second-order convergence is finally obtained along the refinement of time step-sizes for almost all values of $M$.  We also would like to point out that  SSPRK2-SE  with or without SSH reconciliation always fails when $\Delta_{\text{btr}}t$ = 64~s for the tested values of $M>1$. 
When $\Delta_{\text{btr}}t$ = 32~s,  SSPRK2-SE  with SSH reconciliation is also close to explode for all values of $M$ while SSPRK2-SE without SSH reconciliation still works well when $M$ is small. This is probably because that the SSPRK2 substepping for the barotropic mode solve also could not provide good prediction of the SSH  perturbation when the substepping size $\Delta_{\text{btr}}t$ is  large. That is why SSPRK2-SE with SSH reconciliation blows up at $\Delta_{\text{btr}}t$ = 64~s and $M=1$, but SSPRK2-SE without reconciliation does not.

	 The similar behaviors are also observed from Table \ref{baroclinic:m:1:imp} for  SSPRK3-SE when the SSH reconciliation is not present, where the second-order temporal convergence is gradually achieved when $M$ is small  but it  becomes unstable with  deteriorated convergence when $M$ is large.  When the SSH reconciliation is applied, both numerical accuracy and stability of SSPRK3-SE are greatly improved for all cases of $M$. The convergence order along the refinement of time-step sizes is still  around two except when $M=16$ and slightly greater than two when $M=1$. Compared with SSPRK2-SE,  SSPRK3-SE always works very well for $\Delta_{\text{btr}}t$ = 64~s and 32~s, which implies that the SSPRK3 substepping for the barotropic mode solve is more stable and able to produce better prediction of the SSH  perturbation than the SSPRK2 substepping even when the substepping size $\Delta_{\text{btr}}t$ is relatively large and thus the reconciliation process becomes very helpful.	 
 }

		\begin{table}[!ht]
		\centering\footnotesize
		\begin{tabular}{|rc|rrrr|rrrr|}
			\hline
			\multirow{3}{*}{$M$}&\multirow{3}{*}{$\Delta t$-$\Delta_{\text{btr}} t$}& \multicolumn{4}{c|}{Without SSH reconciliation}& \multicolumn{4}{c|}{With SSH reconciliation}\\\cmidrule(lr){3-6}\cmidrule(lr){7-10}
			&&\multicolumn{2}{c}{Velocity}& \multicolumn{2}{c|}{Layer thickness}&\multicolumn{2}{c}{Velocity}& \multicolumn{2}{c|}{Layer thickness}\\ \cmidrule(lr){3-4}\cmidrule(lr){5-6}\cmidrule(lr){7-8}\cmidrule(lr){9-10}
			&&Error &Rate &Error &Rate&Error &Rate &Error &Rate \\ \hline\hline
		\multirow{6}{*}{16}&1024-64&N/A&-&N/A&-&N/A&-&N/A&-\\
		&512-32&N/A&-&N/A&-&2.209E-01& - &1.104E-04& -\\
		&256-16 &N/A&-&N/A&-&1.123E-02&4.30&5.059E-06&4.45\\
		&128-08 &1.486E-02& - &6.613E-06& -&4.886E-03&1.20&2.146E-06&1.24\\
		&64-04 &8.841E-03&0.75&3.686E-06&0.84&1.532E-03&1.67&6.623E-07&1.70\\
		&32-02 &4.435E-03&1.00&1.924E-06&0.94&4.026E-04&1.93&1.745E-07&1.92\\
		\hline
		\multirow{6}{*}{8}&512-64&N/A&-&N/A&-&N/A&-&N/A&-\\
		&256-32&N/A&-&N/A&-&2.194E-01& - &1.097E-04& -\\
		&128-16 &1.463E-02& - &6.685E-06& -&8.653E-03&4.66&3.920E-06&4.81\\
		&64-08 &8.687E-03&0.75&3.629E-06&0.88&2.660E-03&1.70&1.149E-06&1.77\\
		&32-04 &4.344E-03&1.00&1.886E-06&0.94&6.820E-04&1.96&2.957E-07&1.96\\
		&16-02 &1.400E-03&1.63&6.248E-07&1.59&1.714E-04&1.99&7.443E-08&1.99\\
		\hline
		\multirow{6}{*}{4}&256-64&N/A&-&N/A&-&N/A&-&N/A&-\\
		&128-32 &5.769E-01& - &2.796E-03& -&2.381E-01& - &1.260E-04& -\\
		&64-16 &8.081E-03&6.16&3.404E-06&9.68&7.290E-03&5.03&3.265E-06&5.27\\
		&32-08 &3.986E-03&1.02&1.738E-06&0.97&1.824E-03&2.00&7.873E-07&2.05\\
		&16-04 &1.271E-03&1.65&5.676E-07&1.61&4.502E-04&2.02&1.957E-07&2.01\\
		&08-02 &3.304E-04&1.94&1.454E-07&1.96&1.124E-04&2.00&4.888E-08&2.00\\
		\hline
		\multirow{6}{*}{2}&128-64&N/A&-&N/A&-&N/A&-&N/A&-\\
		&64-32 &5.816E-03& - &2.568E-06& -&2.597E-01& - &1.030E-04& -\\
		&32-16 &2.580E-03&1.17&1.153E-06&1.15&6.722E-03&5.27&2.973E-06&5.11\\
		&16-08 &7.568E-04&1.77&3.416E-07&1.76&1.591E-03&2.08&6.878E-07&2.11\\
		&08-04 &1.932E-04&1.97&8.560E-08&2.00&3.911E-04&2.02&1.701E-07&2.02\\
		&04-02 &4.821E-05&2.00&2.117E-08&2.02&9.751E-05&2.00&4.246E-08&2.00\\
		\hline
		\multirow{6}{*}{1}&64-64 &1.163E-02& - &4.929E-06& -&N/A&-&N/A&-\\
		&32-32 &5.415E-03&1.10&2.350E-06&1.07&2.546E-01& - &9.594E-05& -\\
		&16-16 &1.465E-03&1.89&6.282E-07&1.90&6.551E-03&5.28&2.886E-06&5.06\\
		&08-08 &3.632E-04&2.01&1.573E-07&2.00&1.532E-03&2.10&6.624E-07&2.12\\
		&04-04 &9.075E-05&2.00&3.942E-08&2.00&3.761E-04&2.03&1.637E-07&2.02\\
		&02-02 &2.298E-05&1.98&9.970E-09&1.98&9.376E-05&2.00&4.086E-08&2.00\\
		\hline
		\end{tabular}
		\caption{Relative $l_{2}$ errors and convergence rates in $\bu_1$ and $h_1$  for the baroclinic eddies test case produced by the  SSPRK2-SE scheme without vertical mixing.}\label{baroclinic:m:1}
	\end{table}

	\begin{table}[!ht]
		\centering\footnotesize
		\begin{tabular}{|rc|rrrr|rrrr|}
			\hline
			\multirow{3}{*}{$M$}&\multirow{3}{*}{$\Delta t$-$\Delta_{\text{btr}} t$}& \multicolumn{4}{c|}{Without  SSH reconciliation}& \multicolumn{4}{c|}{With SSH reconciliation}\\\cmidrule(lr){3-6}\cmidrule(lr){7-10}
			&&\multicolumn{2}{c}{Velocity}& \multicolumn{2}{c|}{Layer thickness}&\multicolumn{2}{c}{Velocity}& \multicolumn{2}{c|}{Layer thickness}\\ \cmidrule(lr){3-4}\cmidrule(lr){5-6}\cmidrule(lr){7-8}\cmidrule(lr){9-10}
			&&Error &Rate &Error &Rate&Error &Rate &Error &Rate \\ \hline\hline
\multirow{6}{*}{16}&1024-64 &N/A&-&N/A&-&1.251E-02& - &5.173E-06& -\\
&512-32 &N/A&-&N/A&-&1.005E-02&0.32&4.285E-06&0.27\\
&256-16 &N/A&-&N/A&-&7.319E-03&0.46&3.193E-06&0.42\\
&128-08 &1.465E-02& - &6.530E-06& -&4.256E-03&0.78&1.850E-06&0.79\\
&64-04 &8.701E-03&0.75&3.634E-06&0.85&1.913E-03&1.15&8.354E-07&1.15\\
&32-02 &4.354E-03&1.00&1.890E-06&0.94&5.927E-04&1.69&2.566E-07&1.70\\
\hline
\multirow{6}{*}{8}&512-64 &N/A&-&N/A&-&1.128E-02& - &5.048E-06& -\\
&256-32 &N/A&-&N/A&-&7.490E-03&0.59&3.317E-06&0.61\\
&128-16 &1.465E-02& - &6.530E-06& -&4.286E-03&0.81&1.876E-06&0.82\\
&64-08 &8.701E-03&0.75&3.634E-06&0.85&1.919E-03&1.16&8.395E-07&1.16\\
&32-04 &4.354E-03&1.00&1.890E-06&0.94&5.937E-04&1.69&2.571E-07&1.71\\
&16-02 &1.406E-03&1.63&6.273E-07&1.59&1.554E-04&1.93&6.748E-08&1.93\\
\hline
\multirow{6}{*}{4}&256-64 &N/A&-&N/A&-&9.075E-03& - &4.122E-06& -\\
&128-32 &1.464E-02& - &6.541E-06& -&4.785E-03&0.92&2.161E-06&0.93\\
&64-16 &8.695E-03&0.75&3.633E-06&0.85&2.031E-03&1.24&8.984E-07&1.27\\
&32-08 &4.351E-03&1.00&1.888E-06&0.94&6.069E-04&1.74&2.632E-07&1.77\\
&16-04 &1.406E-03&1.63&6.274E-07&1.59&1.565E-04&1.96&6.791E-08&1.95\\
&08-02 &3.665E-04&1.94&1.611E-07&1.96&3.927E-05&1.99&1.708E-08&1.99\\
\hline
\multirow{6}{*}{2}&128-64 &N/A&-&N/A&-&7.200E-03& - &3.313E-06& -\\
&64-32 &8.625E-03& - &3.612E-06& -&3.160E-03&1.19&1.451E-06&1.19\\
&32-16 &4.317E-03&1.00&1.871E-06&0.95&8.300E-04&1.93&3.666E-07&1.98\\
&16-08 &1.405E-03&1.62&6.278E-07&1.58&1.763E-04&2.24&7.665E-08&2.26\\
&08-04 &3.665E-04&1.94&1.613E-07&1.96&4.067E-05&2.12&1.765E-08&2.12\\
&04-02 &9.196E-05&1.99&4.023E-08&2.00&9.890E-06&2.04&4.299E-09&2.04\\
\hline
\multirow{6}{*}{1}&64-64 &8.090E-03& - &3.513E-06& -&6.357E-03& - &2.970E-06& -\\
&32-32 &4.048E-03&1.00&1.757E-06&1.00&2.552E-03&1.32&1.169E-06&1.35\\
&16-16 &1.391E-03&1.54&6.253E-07&1.49&5.055E-04&2.34&2.262E-07&2.37\\
&08-08 &3.662E-04&1.93&1.620E-07&1.95&7.568E-05&2.74&3.358E-08&2.75\\
&04-04 &9.193E-05&1.99&4.031E-08&2.01&1.278E-05&2.57&5.596E-09&2.58\\
&02-02 &2.275E-05&2.01&9.939E-09&2.02&2.635E-06&2.28&1.144E-09&2.29\\
\hline
		\end{tabular}
		\caption{Relative $l_{2}$ errors and convergence rates in $\bu_1$ and $h_1$  for the baroclinic eddies test case produced by the  SSPRK3-SE  scheme without vertical mixing.}\label{baroclinic:m:1:imp}
	\end{table}

\begin{figure}[!ht]
		\centerline{
			\includegraphics[width=1.15in,height=2.3in]{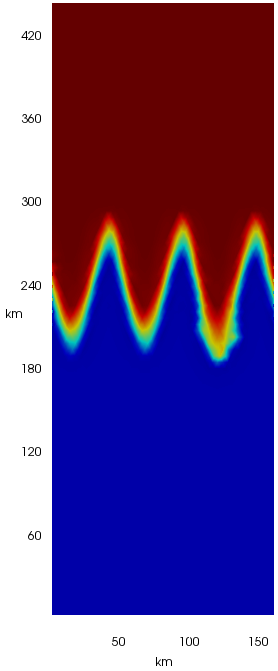}\hspace{-0.1cm}
			\includegraphics[width=1.15in,height=2.3in]{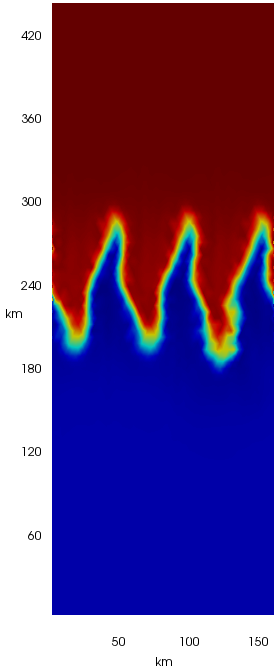}\hspace{-0.1cm}
			\includegraphics[width=1.15in,height=2.3in]{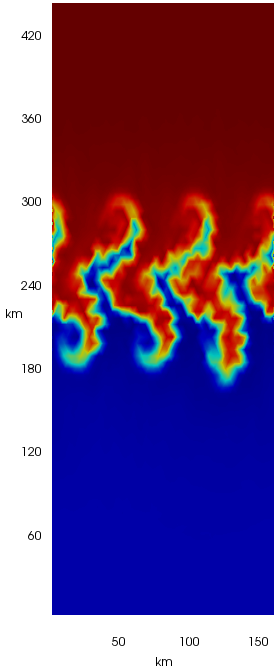}\hspace{-0.1cm}
			\includegraphics[width=1.15in,height=2.3in]{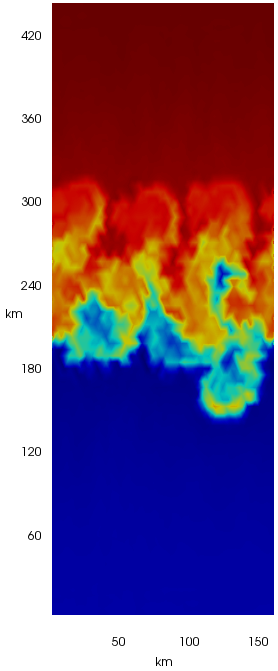}\hspace{-0.1cm}
			\includegraphics[width=1.15in,height=2.3in]{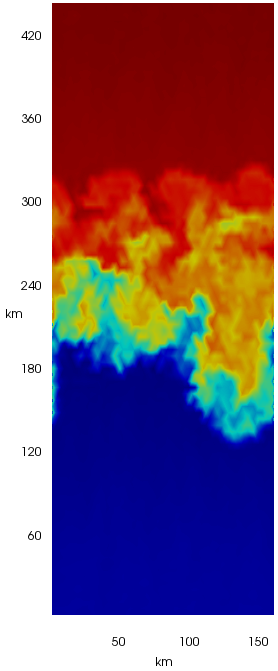}\hspace{-0.1cm}
			\includegraphics[width=0.8in,height=2in]{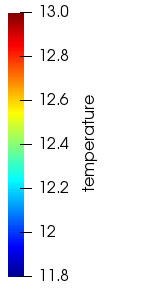}
		} 
		\centerline{
			\includegraphics[width=1.15in,height=2.3in]{./Figures/bot1e-3_day0.png}\hspace{-0.1cm}
			\includegraphics[width=1.15in,height=2.3in]{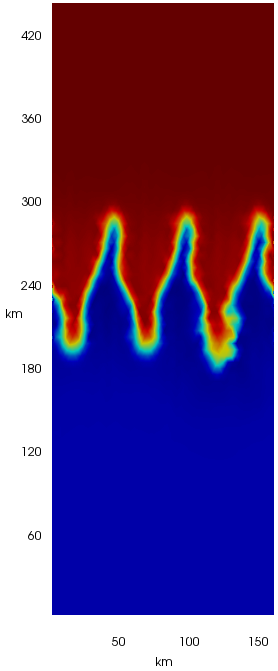}\hspace{-0.1cm}
			\includegraphics[width=1.15in,height=2.3in]{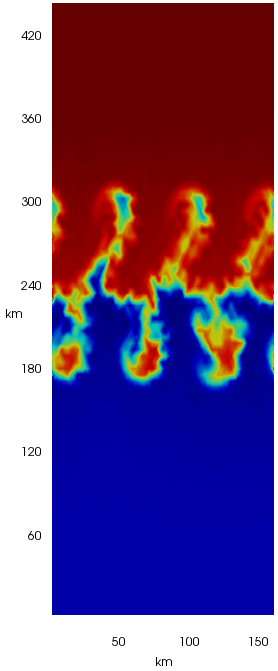}\hspace{-0.1cm}
			\includegraphics[width=1.15in,height=2.3in]{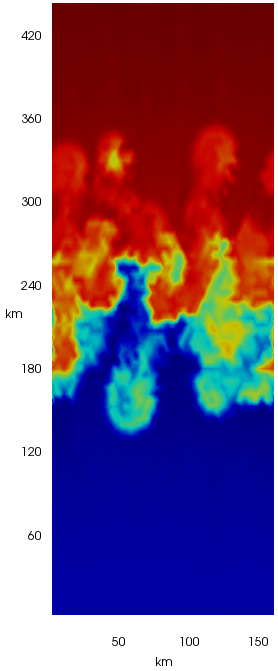}\hspace{-0.1cm}
			\includegraphics[width=1.15in,height=2.3in]{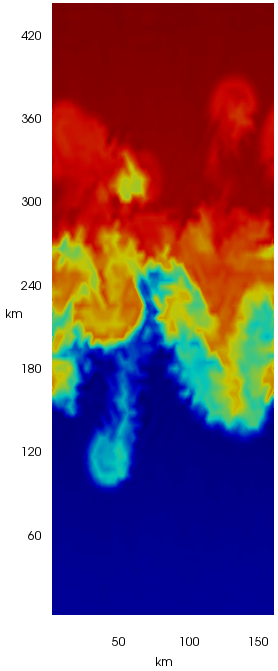}\hspace{-0.1cm}
			\includegraphics[width=0.8in,height=2in]{./Figures/ColorBar_4km.png}
		} 
		\caption{Simulated surface (potential) temperatures by  the SSPRK3-SE  scheme (without vertical mixing) with $(\Delta t, \Delta_{\text{btr}}t)=(240~{\rm s},15~{\rm s})$, i.e., $M=16$ for the baroclinic eddies test case. From left to right are the results at the initial time and Day 5, 10, 20 and 30. Top: $c_{\text{drag}}=0.01$; bottom:  $c_{\text{drag}}=0.001$.}\label{temp:ssprk2}
	\end{figure}

\begin{figure}[!ht]
	\centerline{
		\includegraphics[height=2.5in]{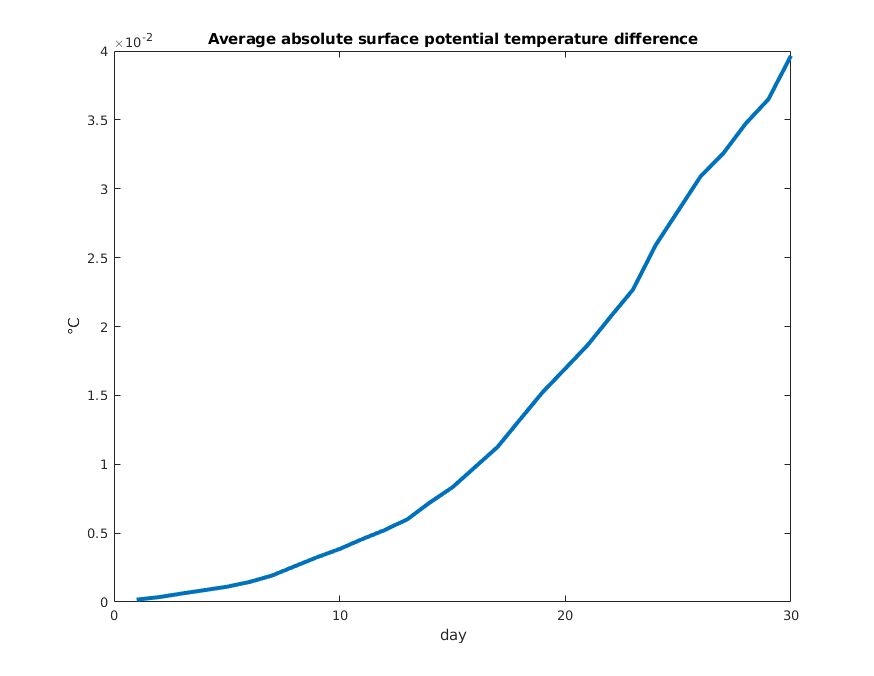}		\hspace{-0.7cm}
		\includegraphics[height=2.5in]{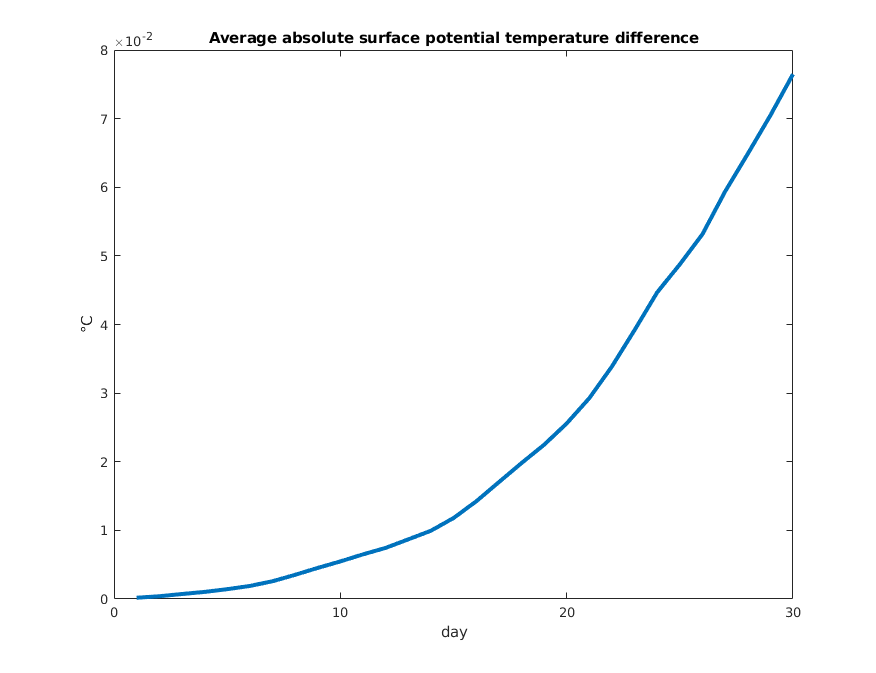}	
} 
	\caption{Evolutions of the average absolute differences in 30 days between the simulated surface temperatures produced  by the SSPRK3-SE  and MPAS-SE schemes (without vertical mixing) with $(\Delta t, \Delta_{\text{btr}}t)=(240~{\rm s},15~{\rm s})$, i.e., $M=16$. Left:  $\ c_{\text{drag}}=0.01$; right: $\ c_{\text{drag}}=0.001$.}\label{temp:diff}
\end{figure}

Based on the  results presented above, we find that the SSH reconciliation process is very important to the proposed SSPRK-SE schemes, so  it will be always applied in all remaining tests.
Finally we take ($\Delta t, \Delta_{\text{btr}} t)$ = (240~s, 15~s), i.e.,  $M=16$ and perform a \textcolor{black}{30-day-long simulation} for the baroclinic eddies test case with different bottom drag coefficients using SSPRK3-SE without vertical mixing. 
In addition, the temperature tracer equation is added by adopting the existing treatment  in MPAS-Ocean. Snapshots of the surface temperature at the initial time and Day 5, 10, 20 and 30 are shown in Figure \ref{temp:ssprk2}, where the top row is associated with the case of $c_{\text{drag}}=0.01$ and the bottom row  with $c_{\text{drag}}=0.001$. It is observed that for the weaker bottom drag, the eddies escape further  away from the interface as is expected. In addition, we also perform the simulation using MPAS-SE and plot the evolutions of the average absolute difference between the simulated surface temperatures produced by  SSPRK3-SE and MPAS-SE in Figure \ref{temp:diff}. 
	We observe that the  differences gradually increase along the time for both  bottom drag choices. Because the eddies travel faster in the smaller bottom drag case, the choice $c_{\text{drag}}=0.001$ results in even 
	larger differences between the results by SSPRK3-SE and MPAS-SE compared to those for the choice $c_{\text{drag}}=0.01$.

\subsection{The global ocean test case with one hundred layers}

	 We next consider a test case having a global real-world ocean configuration from the  MPAS-Ocean platform. The horizontal SCVT mesh, denoted as ``QU240'', is quasi-uniform over the globe, with cell widths of 240~km. There are 100 layers and each of them contains 7,234 cells, 22,736 edges, and 15,459 vertices. As mentioned in Remark \ref{lap-bihar}, the horizontal biharmonic operator with the coefficient $\nu_h=\text{2.0E14}$ is applied as  the hyperviscosity. In this  test case, the bottom drag coefficient  is $c_{\text{drag}}=\text{1.0E-3}$ and the vertical viscosity is $\nu_v=\text{1.0E-4}$. 
	The initial temperature and salinity values at the ocean surface are shown in Figure \ref{QU240:init}.  
	{\color{black}
	We numerically study the temporal accuracy of the proposed schemes by setting the terminal time to
	$T=65536$~s with  $\Delta_{\text{btr}} t =$ 512, 256, 128, 64 and 32~s and 
	$M=$ 32, 16, 8, 4, 2 and 1. The barotropic mode (external gravity waves) is strongly dominating in this case and  the SSH reconciliation  is always applied in the proposed SSPRK-SE schemes. 
The numerical results of SSPRK2-SE and SSPRK3-SE with and without the vertical mixing  are presented in Tables \ref{QU240:m:1:exp} and \ref{QU240:m:1:imp}, respectively. 
\begin{figure}[!ht]
		\centerline{
			\includegraphics[height=1.9in,width=3.6in]{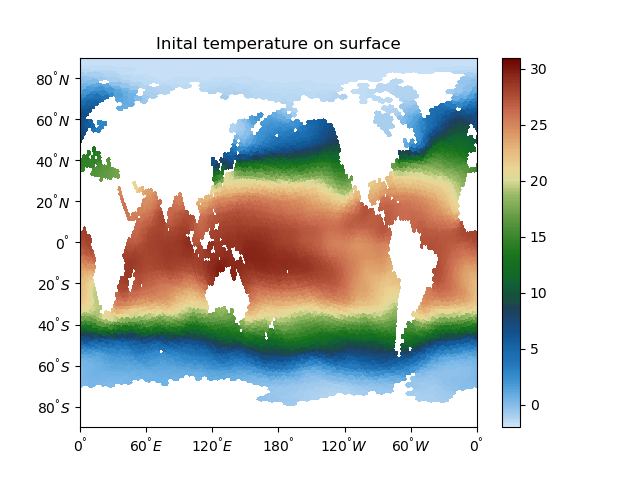}\hspace{-1.2cm}
			\includegraphics[height=1.9in,width=3.6in]{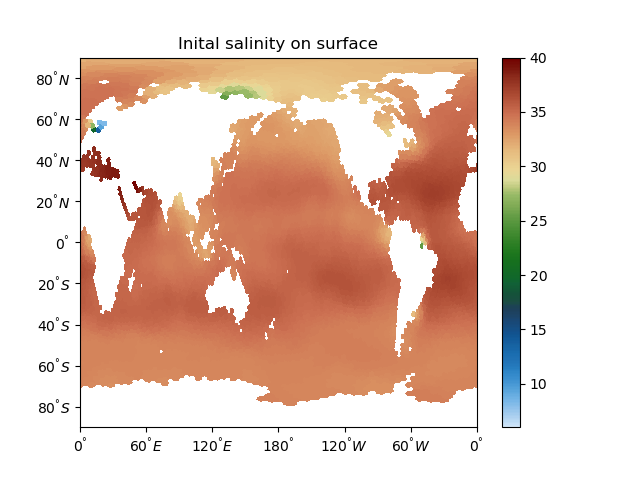}}
			\vspace{-0.2cm}
		\caption{The initial temperature and salinity of the ocean surface for the global ocean test case.}\label{QU240:init}
	\end{figure}

		For SSPRK2-SE, the optimal second-order temporal convergence is obtained along the decrease of the time-step size for both velocity and layer thickness in all the cases when vertical mixing is not applied, but
 the convergence rates of velocity and layer thickness downgrade to the first-order when the vertical mixing is used since the vertical diffusion is quite significant for this test case. We also observe that SSPRK2-SE explodes for all tested values of $M$ when $\Delta_{\text{btr}} t =$ 512~s.
	 The convergence behaviors are a little more complex for SSPRK3-SE. SSPRK3-SE without vertical mixing achieves the second-order  convergence for both velocity and layer thickness. On the other hand, we also observe that the convergence order starts to drop  a little when the relative errors reach 
	around the magnitudes of 1.0E-05 for velocity and and 1.0E-07  for layer thickness. It is probably caused by the tangling effect of the SSH reconciliation process and the spatial discretization error to  the model splitting error when the time-step size is very small but the spatial mesh size is relatively large.
	 With vertical mixing, the convergence rates of velocity and layer thickness again all downgrade to the first-order just as SSPRK2-SE. In addition,
	 the SSPRK3-SE scheme mostly  works very well even when $\Delta_{\text{btr}} t =$ 512~s, which implies that SSPRK3-SE has better numerical stability   than SSPRK2-SE.  The vertical mixing  also  further stabilizes the  SSPRK3-SE scheme  since
	 the simulation blows up without vertical mixing while it does not with vertical mixing when the time-step size is large such as $\Delta_{\text{btr}} t =$ 512~s and $M=32$.
	 Moreover, we observe that SSPRK3-SE  yields  obviously smaller errors than SSPRK2-SE in almost every case  when using the same time-step sizes. All of these again imply the SSPRK3-SE scheme could be a more favorable choice for long time simulations in practice although it is computationally more expensive per time step. }

		\begin{table}[!ht]
		\centering
		\footnotesize\begin{tabular}{|rc|rrrr|rrrr|}
			\hline
			\multirow{3}{*}{$M$}&\multirow{3}{*}{$\Delta t$-$\Delta_{\text{btr}} t$}& \multicolumn{4}{c|}{Without  vertical mixing}& \multicolumn{4}{c|}{With  vertical mixing}\\\cmidrule(lr){3-6}\cmidrule(lr){7-10}
			&&\multicolumn{2}{c}{Velocity}& \multicolumn{2}{c|}{Layer thickness}&\multicolumn{2}{c}{Velocity}& \multicolumn{2}{c|}{Layer thickness}\\ \cmidrule(lr){3-4}\cmidrule(lr){5-6}\cmidrule(lr){7-8}\cmidrule(lr){9-10}
			&&Error &Rate &Error &Rate&Error &Rate &Error &Rate \\ \hline\hline
\multirow{5}{*}{32}&16384-512 &N/A& - &N/A& -&N/A& - &N/A& -\\
&8192-256 &2.690E-01& - &2.397E-04& -&1.843E-01& - &3.933E-04& -\\
&4096-128 &5.434E-02&2.31&5.067E-05&2.24&2.742E-02&2.75&2.045E-04&0.94\\
&2048-64 &1.321E-02&2.04&1.319E-05&1.94&1.568E-02&0.81&1.061E-04&0.95\\
&1024-32 &3.280E-03&2.01&3.478E-06&1.92&1.026E-02&0.61&5.396E-05&0.97\\
\hline
\multirow{5}{*}{16}&8192-512 &N/A& - &N/A& -&N/A& - &N/A& -\\
&4096-256 &5.437E-02& - &5.441E-05& -&2.747E-02& - &2.058E-04& -\\
&2048-128 &1.322E-02&2.04&1.421E-05&1.94&1.569E-02&0.81&1.063E-04&0.95\\
&1024-64 &3.284E-03&2.01&3.732E-06&1.93&1.026E-02&0.61&5.399E-05&0.98\\
&512-32 &8.194E-04&2.00&1.042E-06&1.84&5.794E-03&0.82&2.714E-05&0.99\\
\hline
\multirow{5}{*}{8}&4096-512 &N/A& - &N/A& -&N/A& - &N/A& -\\
&2048-256 &1.330E-02& - &1.915E-05& -&1.575E-02& - &1.073E-04& -\\
&1024-128 &3.303E-03&2.01&4.957E-06&1.95&1.027E-02&0.62&5.415E-05&0.99\\
&512-64 &8.245E-04&2.00&1.326E-06&1.90&5.794E-03&0.83&2.717E-05&0.99\\
&256-32 &2.089E-04&1.98&3.953E-07&1.75&3.047E-03&0.93&1.353E-05&1.01\\
\hline
\multirow{5}{*}{4}&2048-512 &N/A& - &N/A& -&N/A& - &N/A& -\\
&1024-256 &3.513E-03& - &1.108E-05& -&1.033E-02& - &5.525E-05& -\\
&512-128 &8.709E-04&2.01&2.786E-06&1.99&5.799E-03&0.83&2.733E-05&1.02\\
&256-64 &2.201E-04&1.98&7.287E-07&1.93&3.048E-03&0.93&1.356E-05&1.01\\
&128-32 &5.883E-05&1.90&2.084E-07&1.81&1.542E-03&0.98&6.679E-06&1.02\\
\hline
\multirow{5}{*}{2}&1024-512 &N/A& - &N/A& -&N/A& - &N/A& -\\
&512-256 &1.413E-03& - &9.384E-06& -&5.898E-03& - &2.894E-05& -\\
&256-128 &3.359E-04&2.07&2.309E-06&2.02&3.057E-03&0.95&1.379E-05&1.07\\
&128-64 &8.585E-05&1.97&5.846E-07&1.98&1.543E-03&0.99&6.714E-06&1.04\\
&64-32 &2.392E-05&1.84&1.530E-07&1.93&7.568E-04&1.03&3.240E-06&1.05\\
\hline
\multirow{5}{*}{1}&512-512 &N/A& - &N/A& -&N/A& - &N/A& -\\
&256-256 &1.145E-03& - &8.990E-06& -&3.240E-03& - &1.644E-05& -\\
&128-128 &2.609E-04&2.13&2.199E-06&2.03&1.561E-03&1.05&7.087E-06&1.21\\
&64-64 &6.521E-05&2.00&5.496E-07&2.00&7.589E-04&1.04&3.295E-06&1.10\\
&32-32 &1.691E-05&1.95&1.378E-07&2.00&3.558E-04&1.09&1.520E-06&1.12\\
\hline
		\end{tabular}
		\caption{Relative $l_{2}$ errors and convergence rates in $\bu_1$ and $h_1$  for the global ocean test case produced by the proposed SSPRK2-SE scheme with SSH reconciliation.}\label{QU240:m:1:exp}
	\end{table}

	\begin{table}[!ht]
		\centering\footnotesize
		\begin{tabular}{|rc|rrrr|rrrr|}
			\hline
			\multirow{3}{*}{$M$}&\multirow{3}{*}{$\Delta t$-$\Delta_{\text{btr}} t$}& \multicolumn{4}{c|}{Without  vertical mixing}& \multicolumn{4}{c|}{With  vertical mixing}\\\cmidrule(lr){3-6}\cmidrule(lr){7-10}
&&\multicolumn{2}{c}{Velocity}& \multicolumn{2}{c|}{Layer thickness}&\multicolumn{2}{c}{Velocity}& \multicolumn{2}{c|}{Layer thickness}\\ \cmidrule(lr){3-4}\cmidrule(lr){5-6}\cmidrule(lr){7-8}\cmidrule(lr){9-10}
			&&Error &Rate &Error &Rate&Error &Rate &Error &Rate \\ \hline\hline
\multirow{5}{*}{32}&16384-512 &N/A& - &N/A& -&2.753E-01& - &8.619E-04& -\\
&8192-256 &4.590E-02&-&2.257E-04&-&8.758E-02&1.65&4.408E-04&0.97\\
&4096-128 &6.763E-03&2.76&7.417E-05&1.61&4.925E-02&0.83&1.916E-04&1.20\\
&2048-64 &1.336E-03&2.34&1.925E-05&1.95&2.549E-02&0.95&8.250E-05&1.22\\
&1024-32 &3.717E-04&1.85&4.912E-06&1.97&1.280E-02&0.99&3.750E-05&1.14\\
\hline
\multirow{5}{*}{16}&8192-512 &4.590E-02& - &2.257E-04& -&8.759E-02& - &4.408E-04& -\\
&4096-256 &6.764E-03&2.76&7.417E-05&1.61&4.925E-02&0.83&1.916E-04&1.20\\
&2048-128 &1.336E-03&2.34&1.925E-05&1.95&2.549E-02&0.95&8.250E-05&1.22\\
&1024-64 &3.718E-04&1.85&4.911E-06&1.97&1.280E-02&0.99&3.750E-05&1.14\\
&512-32 &1.061E-04&1.81&1.269E-06&1.95&6.375E-03&1.01&1.783E-05&1.07\\
\hline
\multirow{5}{*}{8}&4096-512 &6.782E-03& - &7.418E-05& -&4.925E-02& - &1.916E-04& -\\
&2048-256 &1.341E-03&2.34&1.924E-05&1.95&2.549E-02&0.95&8.250E-05&1.22\\
&1024-128 &3.727E-04&1.85&4.910E-06&1.97&1.280E-02&0.99&3.750E-05&1.14\\
&512-64 &1.062E-04&1.81&1.269E-06&1.95&6.375E-03&1.01&1.783E-05&1.07\\
&256-32 &3.066E-05&1.79&3.521E-07&1.85&3.161E-03&1.01&8.654E-06&1.04\\
\hline
\multirow{5}{*}{4}&2048-512 &1.413E-03& - &1.932E-05& -&2.549E-02& - &8.254E-05& -\\
&1024-256 &3.864E-04&1.87&4.916E-06&1.97&1.280E-02&0.99&3.750E-05&1.14\\
&512-128 &1.078E-04&1.84&1.269E-06&1.95&6.375E-03&1.01&1.783E-05&1.07\\
&256-64 &3.078E-05&1.81&3.519E-07&1.85&3.161E-03&1.01&8.654E-06&1.04\\
&128-32 &1.031E-05&1.58&1.169E-07&1.59&1.555E-03&1.02&4.215E-06&1.04\\
\hline
\multirow{5}{*}{2}&1024-512 &5.519E-04& - &5.240E-06& -&1.280E-02& - &3.756E-05& -\\
&512-256 &1.355E-04&2.03&1.302E-06&2.01&6.375E-03&1.01&1.783E-05&1.07\\
&256-128 &3.310E-05&2.03&3.529E-07&1.88&3.161E-03&1.01&8.655E-06&1.04\\
&128-64 &1.043E-05&1.67&1.167E-07&1.60&1.555E-03&1.02&4.215E-06&1.04\\
&64-32 &4.174E-06&1.32&4.755E-08&1.30&7.522E-04&1.05&2.029E-06&1.05\\
\hline
\multirow{5}{*}{1}&512-512 &3.844E-04& - &2.226E-06& -&6.388E-03& - &1.793E-05& -\\
&256-256 &7.707E-05&2.32&4.580E-07&2.28&3.162E-03&1.01&8.660E-06&1.05\\
&128-128 &1.431E-05&2.43&1.220E-07&1.91&1.555E-03&1.02&4.215E-06&1.04\\
&64-64 &4.363E-06&1.71&4.762E-08&1.36&7.522E-04&1.05&2.029E-06&1.05\\
&32-32 &1.838E-06&1.25&2.106E-08&1.18&3.510E-04&1.10&9.446E-07&1.10\\
\hline
		\end{tabular}
		\caption{Relative $l_{2}$ errors and convergence rates in $\bu_1$ and $h_1$  for the global ocean test case produced by the proposed SSPRK3-SE  scheme with SSH reconciliation.}\label{QU240:m:1:imp}
	\end{table}

\begin{figure}[!ht]	
			\centerline{\includegraphics[height=1.9in]{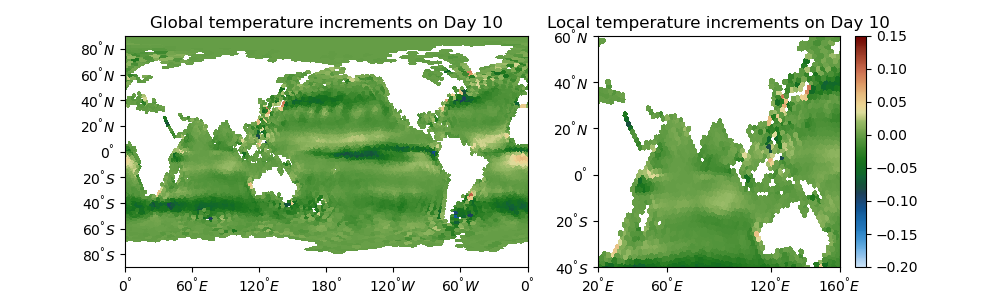}}\vspace{-0.1cm}
			\centerline{\includegraphics[height=1.9in]{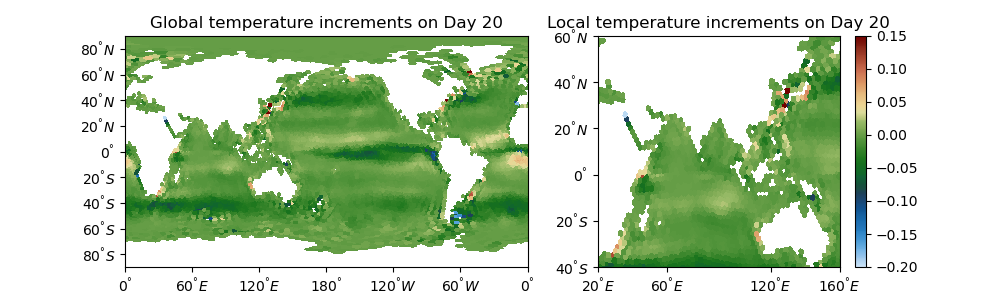}}\vspace{-0.1cm}
			\centerline{\includegraphics[height=1.9in]{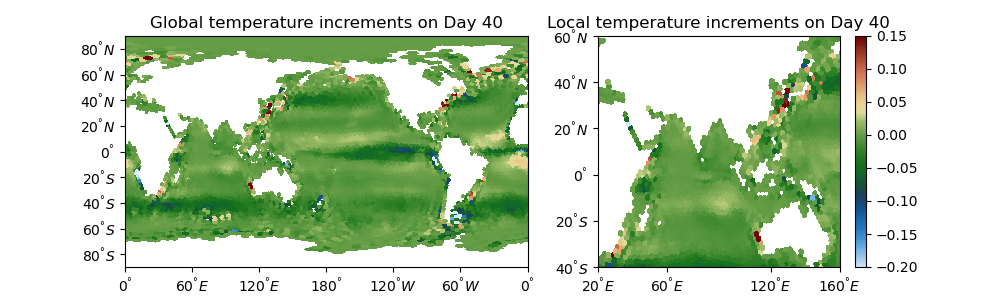}}\vspace{-0.1cm}
			\centerline{\includegraphics[height=1.9in]{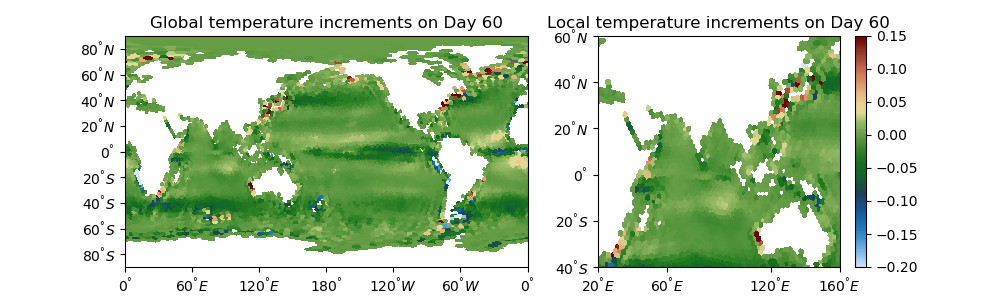}}
		 \vspace{-0.1cm}
		\caption{Simulated surface temperature increments using the SSPRK3-SE scheme (with vertical mixing) with  $(\Delta t, \Delta_{\text{btr}}t)=(3840~{\rm s},240~{\rm s})$, i.e., $M=16$ for the global ocean  test case. From top to bottom are the results on Day 10, 20, 40, and 60. The left ones are the global temperature increments, and the right ones are for the local area at latitude from 40\degree S to 60\degree N, and longitude extending east-ward between 20\degree E and 160\degree E.}\label{temp:ssprk3:QU}
	\end{figure}

\begin{figure}[!ht]
		\centerline{\includegraphics[height=1.9in]{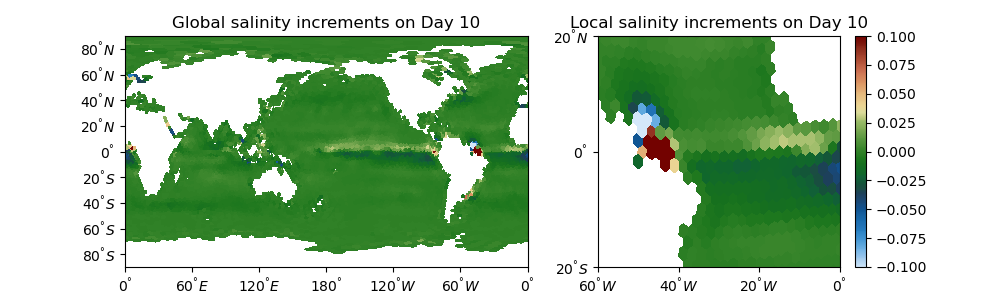}}\vspace{-0.1cm}
			\centerline{\includegraphics[height=1.9in]{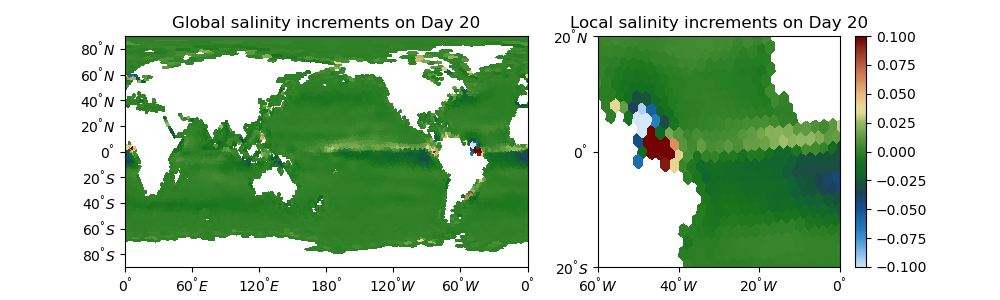}}\vspace{-0.1cm}
			\centerline{\includegraphics[height=1.9in]{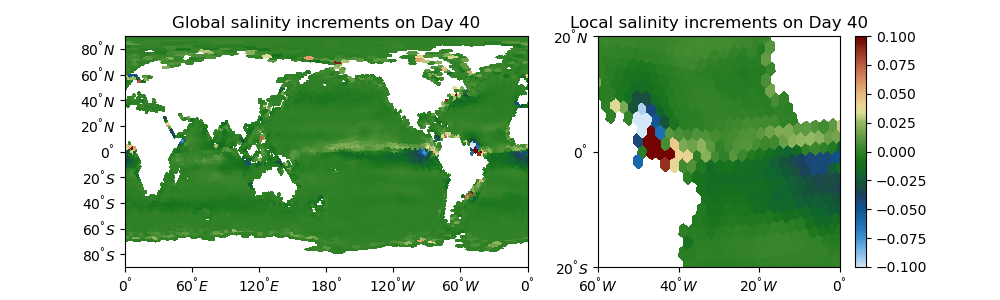}}\vspace{-0.1cm}
			\centerline{\includegraphics[height=1.9in]{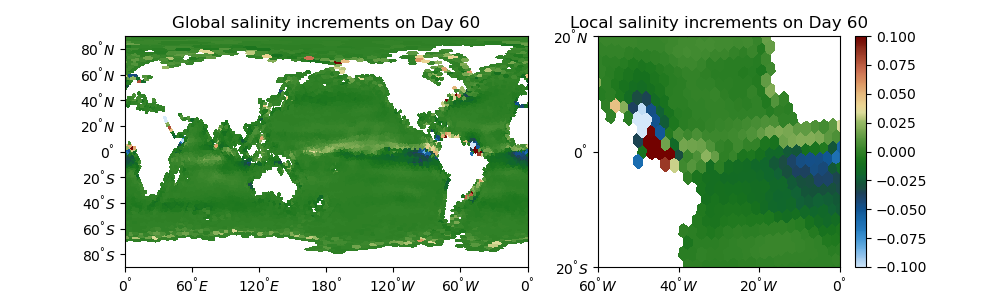}}
		 \vspace{-0.1cm}
		\caption{Simulated surface salinity increments using the SSPRK3-SE scheme (with vertical mixing) with $(\Delta t, \Delta_{\text{btr}}t)=(3840~{\rm s},240~{\rm s})$, i.e., $M=16$ for the global ocean  test case. From top to bottom are the results on Day 10, 20, 40, and 60. The left ones are the global salinity increments, and the right ones are for  the local area at latitude from 20\degree S to 20\degree N, and longitude extending east-ward between 60\degree W and 0\degree E.}\label{sal:ssprk3}
	\end{figure}

\begin{figure}[!ht]
		\centerline{
			\includegraphics[height=2.5in]{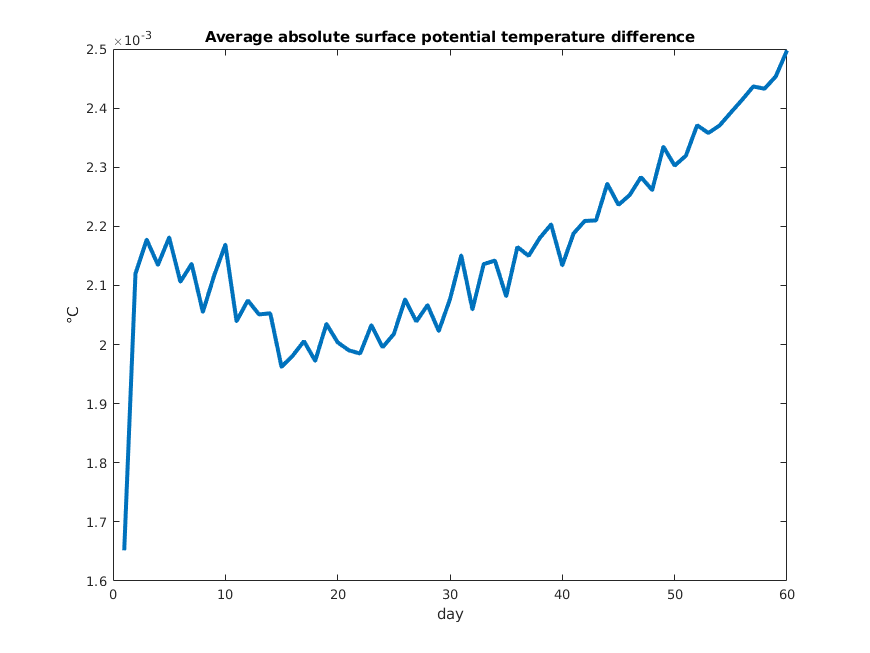}	\hspace{-0.7cm}
			\includegraphics[height=2.5in]{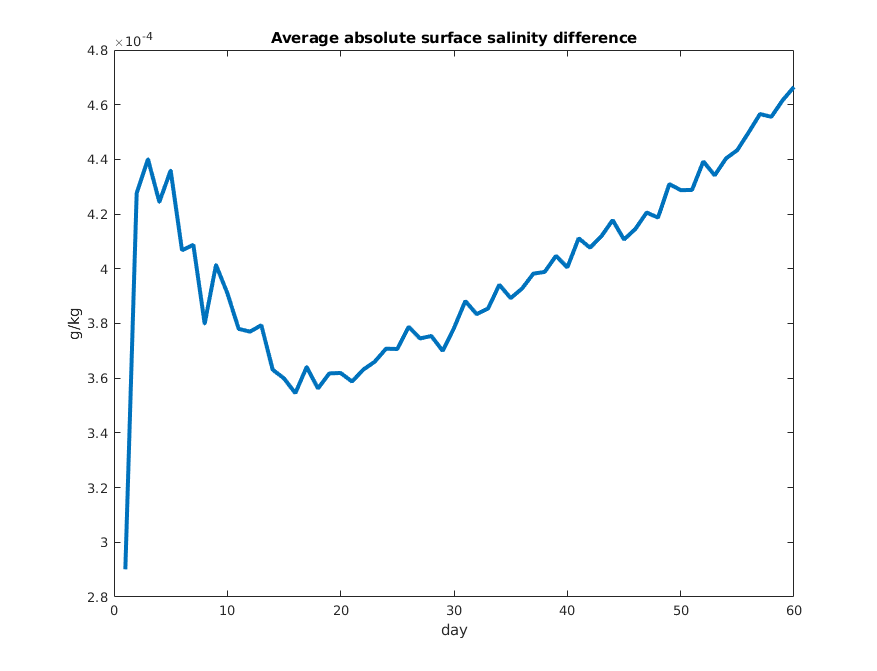}
		} 
		\caption{
		Evolutions of the average absolute differences in 60 days between the simulated surface temperatures (left) and salinity (right) produced by using the SSPRK3-SE and MPAS-SE schemes  (with vertical mixing) with $(\Delta t, \Delta_{\text{btr}}t)=(3840~{\rm s},240~{\rm s})$, i.e., $M=16$.}\label{temp:diff:QU}
	\end{figure}

	By adopting the treatment of tracers from MPAS-Ocean, a \textcolor{black}{60-day-long simulation}  is also carried out for the global ocean test case by using SSPRK3-SE 
	with vertical mixing and ($\Delta t, \Delta_{\text{btr}} t)$ = (3840~s, 240~s), i.e.,  $M=16$. The surface temperature and salinity increments on Day 10, 20, 40, and 60 are shown in Figures \ref{temp:ssprk3:QU} and \ref{sal:ssprk3}. MPAS-Analysis \cite{xylar_asay_davis_2019_2586240} and SciVisColor Colormaps \cite{samsel2data} are used to generate these figures. Besides, we also compare them with the results  produced by MPAS-SE with the same time-step sizes in Figure \ref{temp:diff:QU}. We observe that both of the average absolute differences in  temperature and salinity gradually  increase during the first 10 days, then start to decrease slowly with small oscillations during  the next 10 days, but again keep growing after 20 days until the terminal time.

\subsection{Parallel performance}

	We finally test the parallel performance of the proposed SSPRK2-SE and SSPRK3-SE schemes using the global ocean test case. Following
the MAS-Ocean framework, the domain decomposition is performed only along the horizontal directions, and consequently the unknowns from all layers  but at the same vertical line stay in the same  computing core. In addition to the QU240 mesh, we also consider another finer ocean mesh provided by the MPAS-Ocean platform, denoted there as ``EC60to30'', which contains 60 vertical layers with horizontal resolution varying from 30~km at the equator and poles to 60~km at the mid-latitudes. The grid consists of 235,160 cells, 714,274 edges, and 478,835 vertices at each layer. 

\begin{figure}[!ht]
		\centerline{
			\includegraphics[height=2.4in]{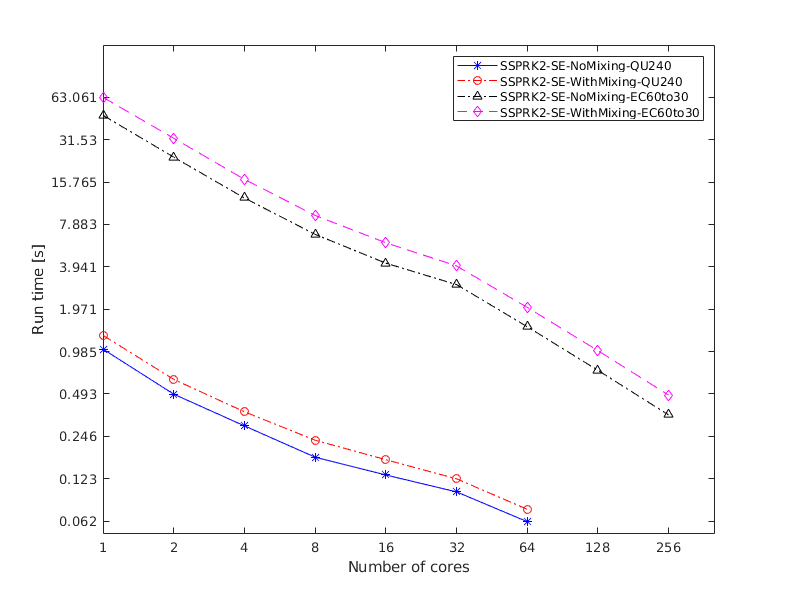}	\hspace{-0.7cm}
			\includegraphics[height=2.4in]{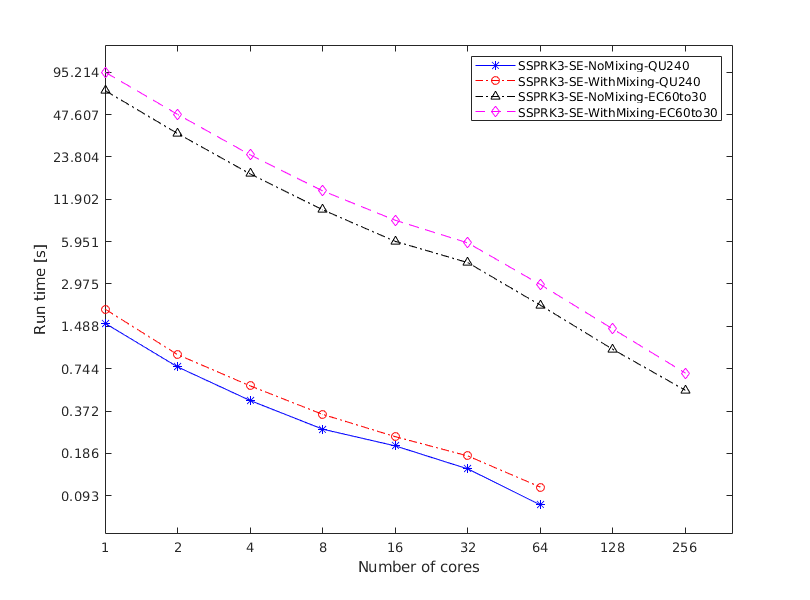}
		} 
		\centerline{
		       \includegraphics[height=2.4in]{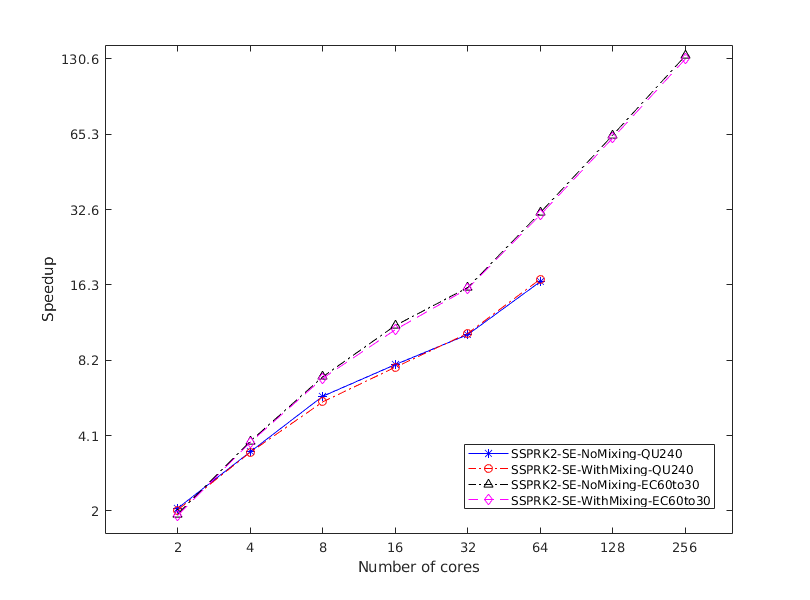}	\hspace{-0.7cm}
			\includegraphics[height=2.4in]{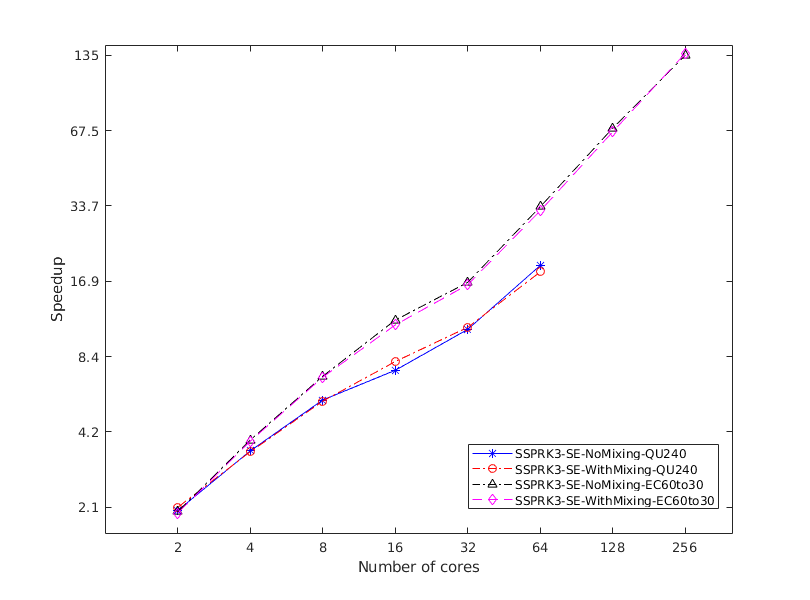}}
		\caption{Plots of the running times per time step vs. the number of cores (top row) and  the corresponding parallel speedups (bottom row) for the global ocean test case using two different global ocean meshes. Left: the SSPRK2-SE scheme with or without mixing; right: the SSPRK3-SE scheme with or without vertical mixing.}\label{parallelperf}
	\end{figure}

We  first fix $\Delta_{\text{btr}}t$ = 15~s and $M=16$  (thus  $\Delta t$ = 240~s) in all simulations.  The tests on the QU240 mesh use up to 64 cores and those on the EC60to30 mesh up to 256 cores. The  plots (in log-log format) of the running time per time step vs. the number of cores and the corresponding speedup plots are presented in Figure \ref{parallelperf}. We see that  in practice the SSPRK-SE schemes with vertical mixing just cost a little more than the SSPRK-SE schemes without vertical mixing. All schemes show very good parallel scalability, i.e., the running times linearly decrease and the speedups linearly increase as the number of cores increases. At the same time, the larger the used spatial mesh is, the better the parallel performance is gained as expected. 

\begin{figure}[!ht]
	\centerline{
		\includegraphics[height=2.4in]{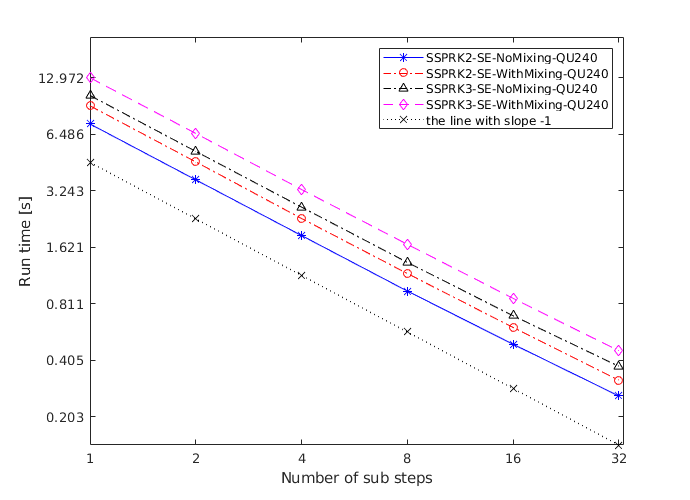}\hspace{-0.7cm}	
		\includegraphics[height=2.4in]{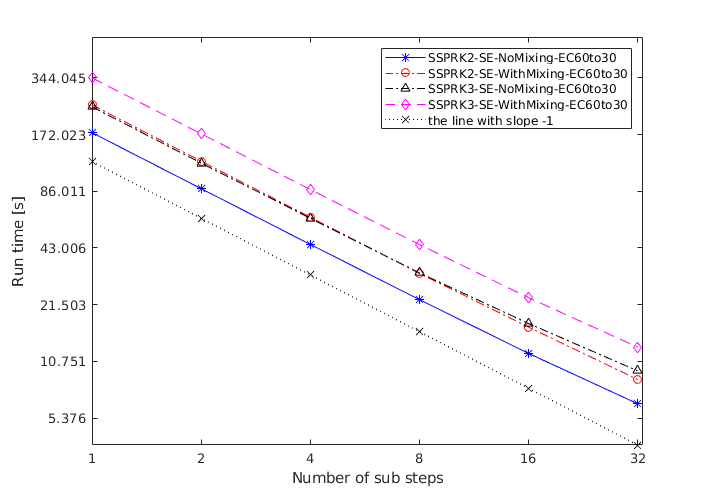}}
	\caption{Plots of the running times vs. the number of substepps for the global ocean test case using two different global ocean meshes. Left: QU240 mesh; right: EC60to30 mesh.}\label{differentMs}
\end{figure}

{\color{black}
 Next we investigate the computational cost and scalability of the proposed SSPRK-SE schemes under different numbers of substeps $M$ for the barotropic mode solve.  Specifically, we  fix the terminal time $T=1920$~s and  $\Delta_{\text{btr}}t$ =15~s, then take  $M=$ 1, 2, 4, 8, 16 and 32 (thus  $\Delta t= 15M$~s) to run all the tests with 64 cores.
 Since the cost of the baroclinic mode solves dominates that of the barotropic mode solves  for this problem when $L>>M$, it is expected that
 the total cost is approximately decreased half when $M$ is doubled and not very large.  The plots (in log-log format) of the running time vs. the number of substeps in all cases are presented in Figure \ref{differentMs}. It is observed that all the curves do roughly  follow the inverse linear relation (i.e., with a slope of $-1$), which demonstrates the benefit and efficiency of the proposed multirate explicit  time-stepping schemes over the global uniform 
 time-stepping.}
	
	\section{Conclusions}\label{sect:conclusion}

	This paper is concerned with high-order (greater than one) numerical methods for  the barotropic  and baroclinic dynamic split system for the layered primitive equations.  Two multirate explicit time-stepping schemes (SSPRK2-SE and SSPRK3-SE) are designed and analyzed based on the framework of classic SSPRK approach, in which  a large time step can be used for the three-dimensional baroclinic  mode solve and a small time step for the two-dimensional barotropic mode solve, and furthermore, each of the two mode solves just need satisfy their respective CFL conditions for numerical stability. Extensive numerical tests on two benchmark test problems from  the MPAS-Ocean platform are performed to demonstrate their high-order accuracy, stability and parallel scalability.  On the other hand, there are still some important questions and tasks worthy  of further  study. \textcolor{black}{
A thorough  convergence analysis of the proposed SSPRK-SE schemes in the fully (both time and space) discrete settings is still open and
but quite sophisticated  since the errors from time integration, spatial discretization, mode splitting and even vertical mixing are tangled up together. Furthermore, how to modify the proposed SSPRK3-SE   to recover the third-order accuracy represents another interesting and important question.}
 It also would be highly useful in practice to develop
appropriate operator splitting  techniques for the vertical mixing  so that  the overall high-order accuracy is not affected when combined with the proposed SSPRK-SE schemes. In the end,  incorporating tracer solves into the proposed schemes without losing high-order temporal accuracy also remains an important task to consider in the future.
	
	\section*{Acknowledgement}
	We sincerely thank the anonymous referees very much for their support and insightful comments which greatly improved the paper. This work was supported by the U.S. Department of Energy, Office of Science, Office of Biological and Environmental Research through Earth and Environmental System Modeling and Scientific Discovery through Advanced Computing programs at Los Alamos National Laboratory and under university grants DE-SC0020270 and DE-SC0020418. This research used the computing resources of National Energy Research Scientific Computing Center, a U.S. Department of Energy Office of Science User Facility operated under Contract No. DE-AC02-05CH11231.
	
	\bibliographystyle{elsarticle-num}
	\bibliography{ssprkRef}
\end{document}